\documentclass[aap,MSNbibl,citesort,dvips]{arximspdf}

\doi{10.1214/09-AAP672}
\volume{20}
\issue{5}
\pubyear{2010}
\firstpage{1761}
\lastpage{1800}

\makeatletter

\newproclaim{definicion}{Definition}[section]
\newproclaim{rem}[definicion]{Remark}
\newtheorem{lema}[definicion]{Lemma}
\newtheorem{teorema}[definicion]{Theorem}
\newtheorem{proposicion}[definicion]{Proposition}
\newtheorem{corolario}[definicion]{Corollary}

\def\tn{\Vert\hspace*{-1.4pt}\vert}

\makeatother

\begin{document}
\begin{frontmatter}

\title{Stochastic vortex method for forced three-dimensional
Navier--Stokes equations and~pathwise convergence rate\thanksref{T2}}
\runtitle{3d-vortex method for forced Navier--Stokes equations}

\thankstext{T2}{Supported by Fondecyt 1070743 and
FONDAP-BASAL-CONICYT.}

\begin{aug}
\author[A]{\fnms{J.} \snm{Fontbona}\corref{}\ead[label=e1]{fontbona@dim.uchile.cl}}
\runauthor{J. Fontbona}
\affiliation{Universidad de Chile}
\address[A]{Departamento de Ingenier\'{\i}a Matem\'{a}tica\\
\quad y Centro de Modelamiento Matem\'{a}tico\\
UMI(2807) UCHILE-CNRS\\
FCFM, Universidad de Chile\\
Casilla 170-3, Correo 3\\
Santiago\\ Chile\\
\printead{e1}} 
\end{aug}

\received{\smonth{10} \syear{2008}}
\revised{\smonth{12} \syear{2009}}

%
\begin{abstract}
We develop a McKean--Vlasov interpretation of Navier--Stokes equations
with external force field in the whole space, by associating with local
mild $L^p$-solutions of the 3d-vortex equation a generalized nonlinear
diffusion with random space--time birth that probabilistically describes
creation of rotation in the fluid due to nonconservativeness of the
force. We establish a local well-posedness result for this process and
a stochastic representation formula for the vorticity in terms of a
vector-weighted version of its law after its birth instant. Then we
introduce a stochastic system of 3d vortices with mollified interaction
and random space--time births, and prove the propagation of chaos
property, with the nonlinear process as limit, at an explicit pathwise
convergence rate. Convergence rates for stochastic approximation
schemes of the velocity and the vorticity fields are also obtained. We
thus extend and refine previous results on the probabilistic
interpretation and stochastic approximation methods for the nonforced
equation, generalizing also a recently introduced random
space--time-birth particle method for the 2d-Navier--Stokes equation
with force.
\end{abstract}

%
\begin{keyword}[class=AMS]
\kwd[Primary ]{60K35}
\kwd{65C35}
\kwd{76M23}
\kwd{76D17}
\kwd[; secondary ]{35Q30}.
\end{keyword}
\begin{keyword}
\kwd{3d-Navier--Stokes equation with external force}
\kwd{McKean--Vlasov model with random space--time birth}
\kwd{stochastic vortex method}
\kwd{propagation of chaos}
\kwd{convergence rate}.
\end{keyword}

\end{frontmatter}

\section{Introduction}

The Navier--Stokes equation for a homogeneous and incompressible
fluid in the whole plane or space, subject to an external force
field $\mathbf{F}$, is given by
%
%
\begin{eqnarray}\label{NSF}
\frac{{\partial}\mathbf{u}}{{\partial t}} + (\mathbf{u}\cdot
\nabla
)\mathbf{u}&=&\nu\Delta\mathbf{u}-\nabla\mathbf{p}+\mathbf{F};\nonumber\\[-8pt]\\[-8pt]
\operatorname{div} \mathbf{u}(t,x)&=&0; \qquad\mathbf{u}(t,x)\to0 \qquad\mbox{as }
|x|\to
\infty.\nonumber
\end{eqnarray}
Here, $\mathbf{u}$ denotes the velocity field, $\mathbf{p}$ is the (unknown)
pressure function and $\nu>0$ is the (constant) viscosity
coefficient. When $\mathbf{F}=0$ or, more generally, when $\mathbf
{F}=\nabla\Psi$
is a conservative field, a probabilistic interpretation
of (\ref{NSF}) in space dimension two
was first developed in 1982 by Marchioro and Pulvirenti \cite{MaPu}.
Their approach was based on
the vortex equation satisfied by the (scalar) field $\operatorname{curl} \mathbf{u}$,
which in 2d and for the case of a conservative external field,
was interpreted as a nonlinear Fokker--Planck (or McKean--Vlasov)
equation with signed initial condition. This was associated with
a nonlinear diffusion process in the sense of McKean, involving
singular interactions through the kernel of Biot--Savart. (For a
general background on the McKean--Vlasov model, we refer the reader
to Sznitman \cite{Szn} and M\'{e}l\'{e}ard \cite{Mel0}.) This approach
led them to the definition of a stochastic system of particle or
vortices with ``mollified'' mean field interaction, for which the
time-marginal empirical measures converge to a solution of the
vortex equation associated with (\ref{NSF}). The convergence on
the path space of that particles system (or, equivalently, the
propagation of chaos property) was proved later by M\'{e}l\'{e}ard
in~\cite{Mel1}. Those works provided a rigorous mathematical meaning
of Chorin's vortex algorithm, heuristically proposed in
\cite{Cho0} as a probabilistic method to simulate the solution of
the 2d-Navier--Stokes equation (see also \cite{Cho}).

In dimension 3, the vorticity field $\mathbf{w}=\operatorname{curl} \mathbf{u}$ is a
solution of
the vectorial nonlinear equation
%
%
\begin{eqnarray}\label{3DVF}
\frac{{\partial}\mathbf{w}}{{\partial t}} + (\mathbf{u}\cdot
\nabla
)\mathbf{w}&=&(\mathbf{w}\cdot\nabla)\mathbf{u}+\nu\Delta\mathbf
{w}+\mathbf{g},\nonumber\\[-8pt]\\[-8pt]
\operatorname{div} w_0&=&0,\nonumber
\end{eqnarray}
where $\mathbf{g}=\operatorname{curl}\mathbf{F}$ and where the relation
%
%
\begin{equation}\label{3BSlaw}
\mathbf{u}(t,x)=\mathbf{K}(\mathbf{w})(t,x):=-\frac{1} {4\pi}\int
_{\mathbb{R}^3}
\frac{(x-y)}{|x-y|^3}\wedge\mathbf{w}(t,y) \,dy
\end{equation}
holds, thanks to the
incompressibility condition $\operatorname{div}\mathbf{u}=0$ and the Biot and Savart law.
Here, $\wedge$ stands for the vectorial
product in $\mathbb{R}^3$, $K(x)\wedge:=-\frac{1}{4\pi}\frac
{x}{|x|^3}\wedge$ is the
three-dimensional Biot--Savart kernel and $\mathbf{K}$ is the Biot--Savart
operator in 3d. (We refer to Bertozzi and Majda \cite{BerMa} for
this and for background on vorticity.)

In absence of external forces, the problem of proving the
approximation of solutions of the 3d-Navier--Stokes equations by a
stochastic system of mean field interacting particles was first
addressed by Esposito and Pulvirenti \cite{EsPu}. In that work, an
approximation result of local solutions by a stochastic system
of three-dimensional vortices with cutoff and mollified
interactions was obtained for each time instant, for initial
vorticities that belonged to $L^1$ together with their Fourier
transform. The convergence held for mollifying parameters that
depended on the realizations of the empirical measures of the
paths of the driving Brownian motions.

Recently, we considered in \cite{F1} the mild version of the
3d-vortex equation with $\mathbf{g}=0$ in the $L^p$ spaces for $p
> \frac{3}{2}$. We proved local (in time) well posedness and regularity
results for that equation, and, under an additional $L^1$
assumption on $w_0$, we showed the equivalence between such
solutions and a generalized nonlinear McKean--Vlasov process with values
in $\mathbb{R}^3\times\mathbb{R}^{3\otimes
3}$ and singular drift term at $t=0$.
We then introduced a system of stochastic 3d vortices with cutoff and
mollified interaction, and
proved the pathwise propagation of chaos property with as limit
the nonlinear process,
deducing moreover stochastic particle approximation results for
the velocity and vorticity fields. (We refer to \cite{FAP} for a
rectification of the discussion in \cite{F1} about the work
\cite{EsPu}.) During the preparation of this work, we have also
become aware of the more recent work of Philipowski \cite{Phi},
who obtained (also in the case $\mathbf{g}=0$) a convergence rate for a
mean field particle approximation of the vorticity field, for a
simpler variation of the system introduced in \cite{F1}. (The
pathwise propagation of chaos property was not addressed.)

In presence of an external force field, the
additional additive term $\mathbf{g}=\operatorname{curl}\mathbf{F}$ in the (2d or
3d) vortex
equation is
physically interpreted as creation of rotation in the fluid. In order to
describe this phenomenon probabilistically, a nonlinear
McKean--Vlasov diffusion process with random space--time
birth was recently
associated with the 2d-vortex equation in Fontbona and M\'{e}l\'{e}ard
\cite{FM}. More precisely, the
law $P_0(dt,dx)$ of the instant and position of birth was suitable, defined
in terms of the initial vorticity and of the external field
$\operatorname{curl}\mathbf{F}$, and it was shown that a scalar-weighted version of the
time marginal law of this process after its birth time was equal
to the solution to the 2d-vortex equation (with $L^1$ data) in a
given interval. The propagation of chaos property was established
for an approximating system of interacting vortices, which were
given birth independently at random positions and times following
the law $P_0$, and a pathwise convergence rate was obtained under
slight additional integrability assumptions on the data.

The first purpose of the present paper is to extend the results
of \cite{F1} and \cite{FM} to the 3d-Navier--Stokes equation with
nonconservative external force field. More precisely, fix $T>0$
and assume that $w_0\dvtx\mathbb{R}^3\to\mathbb{R}^3$ and $\mathbf
{g}\dvtx\mathbb{R}^3\times[0,T]
\to\mathbb{R}^3$ are divergence-free $L^1$-fields. Denote by $I_3$ the
identity matrix in $\mathbb{R}^3$ and let $(B_t)$ be a standard
3d-Brownian motion. Our main goal will be to study the well posedness
on $[0,T]$ of the following nonlinear process, with singular
interaction kernel and values in $\mathbb{R}^3\times\mathbb
{R}^{3\otimes3}$:
%
%
\begin{eqnarray}\label{3elprocesointro}
X_t &=& X_0+\sqrt{2\nu} \int_0^t \mathbf{1}_{\{s\geq\tau\}}\, dB_s
+\int_0^t
\mathbf{K}(\tilde{\rho})
(s,X_s)\mathbf{1}_{\{s\geq\tau\}}\,ds ,\nonumber\\[-8pt]\\[-8pt]
\Phi_t &=& I_3+\int_0^t \nabla\mathbf{K}(\tilde{\rho})(s,X_s)\Phi
_s \mathbf{1}_{\{
s\geq\tau\}} \,ds,\nonumber
\end{eqnarray}
where: $(\tau,X_0)$ is a
random variable in $[0,T]\times\mathbb{R}^3$ (independent of $B$) with
law
\[
P_0(dt,dx)\propto\delta_0(dt) |w_0(x)| \,dx + |\mathbf{g}(t,x) |\,dx\,dt,
\]
$\tilde{\rho}=\tilde{\rho}(t,x)$ is defined for each $t$ from the
law of $(\tau,X,\Phi)$ as
%
%
\begin{equation}\label{lequiv}
\int_{\mathbb{R}^3}\mathbf{f}(y)\tilde{\rho}(t,y)\,dy :=E\bigl(\mathbf{f}(X_t)
\Phi_t
h(\tau,X_0)\mathbf{1}_{\{t\geq\tau\}}\bigr)\qquad \mbox{for } \mathbf
{f}\dvtx\mathbb{R}^3\to
\mathbb{R}^3,
\end{equation}
and $h$ in (\ref{lequiv}) is the density with respect to $P_0$ of
the vectorial measure
$\delta_0(dt)\times\break w_0(x)\,dx + \mathbf{g}(t,x)\,dx\,dt$. [We observe that it is
(\ref{3elprocesointro}) \textit{together} with relation
(\ref{lequiv}) that specify a ``nonlinear process'' in McKean's
sense.]

As we shall see, there will exist a correspondence between mild
$L^p(\mathbb{R}^3) \cap L^1(\mathbb{R}^3)$-solutions $\mathbf{w}$ of
(\ref{3DVF}) for
$p>\frac{3}{2}$, and suitable solutions of the nonlinear
stochastic differential equation
(\ref{3elprocesointro}) and (\ref{lequiv}), through the relation
$\mathbf{w}=\tilde{\rho}$. Thus,
(\ref{lequiv}) provides a representation formula for solutions
$\mathbf{w}$ of (\ref{3DVF}) which extends the one obtained in \cite{F1}
when $\mathbf{g}\equiv0$ (or $\tau\equiv0$). In the present case, this
representation can be intuitively understood as follows. A point
vortex is given birth at random instant and position $(\tau,X_0)$,
rotating in direction $h(\tau,X_0)\in\mathbb{R}^3$. It then evolves
under the effect of diffusion and of the velocity field $\mathbf
{K}(\mathbf{w})$
in (\ref{3elprocesointro}), while its rotation direction and
magnitude are changed under the action of the matrix process
$\Phi_t$ which accounts for the vortex stretching proper to
dimension $3$. Averaging the rotation vectors on the position of
infinitely ``already born vortices'' yields a macroscopic
vorticity field $\mathbf{w}(t)=\tilde{\rho}(t)$, weakly defined by
(\ref{lequiv}). The velocity field instantaneously experienced by
each individual vortex is finally recovered from
$\mathbf{w}$ as a mean field effect through the
interaction kernel of Biot--Savart.

We will adapt the ideas and analytic techniques in \cite{F1} to
first establish local well-posedness and regularity results for
the mild formulation of the vortex equation. Based on this, we
shall then prove local [i.e., for small enough $T>0$ or data
$(w_0,\mathbf{g})$] pathwise well posedness for the nonlinear stochastic
differential equation
(\ref{3elprocesointro}) and (\ref{lequiv}), which will have singular
drift terms at $t=0$.

We shall then introduce a stochastic system of $n$ particles in
$\mathbb{R}^3\times\mathbb{R}^{3\otimes3}$ (or 3d-vortices) with
cutoff and
mollified interaction kernels, and with random space--time births.
The second goal of this paper will be to prove the
strong pathwise convergence of each of these particles as $n$
goes to $\infty$, towards the nonlinear process, at an explicit
rate. To that end, we will improve the techniques used in
\cite{F1} to study the nonlinear process, which relied on
tightness estimates for approximating processes and martingale
problem characterization. More precisely, by a fine use of
regularity properties of the equation, and inspired by ideas
introduced in \cite{FM}, we will show that the approximating
``mollified processes''
converge pathwise at the same rate at which mollified
versions of the vortex equation converge to the original one. We
will be able to exhibit that rate for a large class of
mollified kernels, thanks to classic regularization
techniques in Raviart \cite{Rav} (which are also similar to those
used in \cite{Phi}). These results will imply the propagation of
chaos in a strong norm and, classically, an explicit rate in some
pathwise Wasserstein distance $\mathcal{W}$. From this we will also
deduce convergence rates for approximation schemes of the
vorticity and velocity fields. Unfortunately, the mollifying
parameter will be required to go very slowly to $0$ as $n$ goes to
$\infty$, which will yield a very slow (but not necessarily
optimal) rate for the particles convergence.

Finally, we point out that our regularity results on the mild
equation in $L^p$ will ensure that the stochastic flow
%
%
\begin{equation}\label{3linearflow}
\xi_{s,t}(x)=x+ \sqrt{2\nu}(B_t-B_s)+
\int_s^t\mathbf{u}(r,\xi_{s,r}(x))\,dr
\end{equation}
is of class $C^1(\mathbb{R}^3)$, and so one can write
%
%
\begin{equation}\label{formstocflow}
(X_t,\Phi_t)\mathbf{1}_{\{t\geq\tau\}}=(\xi_{\tau,t}(X_0),\nabla_x
\xi_{\tau,t}(X_0))\mathbf{1}_{\{t\geq\tau\}}.
\end{equation}
Equation
(\ref{lequiv}) can thus be seen as a stochastic analog for the
3d-Navier--Stokes equation of the ``Lagrangian representation'' of the
vorticity of the 3d-Euler equation $\nu=0$ (see, e.g., \cite{ChoMa},
Chapter 1), an analogy established in \cite{EsPu,F1} when
$\mathbf{g}
\equiv0$. Lagrangian representations of the 3d-Navier--Stokes
equations as stochastic analogues to representations formulae for
the Euler equation have been studied by several authors, some of
which have led to (local) well-posedness results for the equation.
See, for example, Esposito et al. \cite{EsMaPuSc} and, for more recent
developments,
Busnello et al. \cite{BuFraRo} and Iyer \cite{I}. The latter
works follow approaches that are in some sense ``dual'' to ours,
establishing representations of strong solutions of the vortex or
Navier--Stokes equations in terms of expectations of the initial
data, after being transported and modified by the stochastic
flow. A related stochastic approach is adopted in Gomes \cite{Go}
to establish a variational formulation of the Navier--Stokes
equation, analogous to Arnold's variational characterization of
the Euler equation. A seemingly very different further
probabilistic point of view, providing global well posedness for
small initial data, was introduced by Le Jan and Sznitman in
\cite{LJS}, who associated with the Fourier transform of the
velocity field a multitype branching process or stochastic
cascade. See, for example, Bhattacharya et al. \cite{Bha} for more recent
developments in that direction.

The remainder of this work is organized as follows. In Section \ref{sec2}
we first present a weak formulation of
(\ref{3elprocesointro}) and (\ref{lequiv}) in terms of a nonlinear
martingale problem, and discuss its connection with
(\ref{3DVF}). In Section \ref{sec3}, we shall obtain local well-posednes
and regularity results for the mild version of the vortex
equation in $L^p$, for $p \in(\frac{3}{2},3)$. In Section \ref{sec4} we
state some results about a nonlinear Fokker--Planck equation with
external field associated with the process with random space--time
birth $X$ in (\ref{3elprocesointro}). We use this and the
previous results to show strong local-in-time well posedness for
the nonlinear stochastic differential equation
(\ref{3elprocesointro}) and (\ref{lequiv}). We, moreover, obtain the
pathwise convergence result and estimates for approximating
mollified versions of that problem. In Section \ref{sec5}, we introduce the
system of 3d-stochastic vortices with random space--time birth, and
deduce the propagation of chaos property and its rate. We also
prove approximation results for the velocity and the vorticity of
the forced 3d-Navier--Stokes equation with their corresponding
convergence rates. In Section \ref{sec6} we shall discuss how these
rates of convergence are slightly improved when Sobolev
regularity of the initial condition and external field is assumed.

Let us establish some notation:

\begin{enumerate}[--]
\item[--] By $\mathcal{M}\mathit{eas}^T$ we denote
the space of measurable real-valued
functions on $[0,T]\times\mathbb{R}^3$.

\item[--] $C^{1,2}$
is the set of real-valued
functions on $[0,T]\times\mathbb{R}^3$ with
continuous derivatives up to the first order in $t\in[0,T]$ and up
to the second order in $x\in\mathbb{R}^3 $. $C_b^{1,2}$
is the subspace of bounded functions in $C^{1,2}$
with bounded derivatives.

\item[--] $\mathcal{D}$ is the space of compactly supported functions
on $\mathbb{R}^3$ with
infinitely many derivatives.

\item[--] For all $1\leq p \leq\infty$ we denote by
$L^p$ the space $L^p(\mathbb{R}^3)$ of real-valued functions on
$\mathbb{R}^3$.
By $\|\cdot\|_p$ we denote the corresponding norm, and $p^*$ stands
for the H\"{o}lder conjugate of $p$. We write
$W^{1,p}=W^{1,p}(\mathbb{R}^3)$ for the Sobolev space of functions in
$L^p$ with partial derivatives of first order in $L^p$.

\item[--] If $E$ is a space of real-valued functions (defined on
$\mathbb{R}^3$ or on $[0,T]\times\mathbb{R}^3$), then the notation
$(E)^3$ is
used for the space of $\mathbb{R}^3$-valued functions with scalar
components in $E$. If $E$ has a norm, the norm in $(E)^3$ is
denoted in the same way.

\item[--] For notational simplicity, if $\mathbf{f},\mathbf
{g}\dvtx\mathbb{R}^3\to\mathbb{R}
^3$ are vector
fields and $Z\dvtx\mathbb{R}^3\to\mathbb{R}^{3\otimes3}$ is a matrix
function, we
will write $\mathbf{f}\mathbf{g}:=\sum_i^3\mathbf{f}_i\mathbf
{g}_i$ and $\mathbf
{f}Z$ for the row-vector
$(\mathbf{f}^tZ)_i:= \sum_{j=1}^3 \mathbf{f}_j Z_{j,i}$. By $\nabla
\mathbf{f}$ we denote\vspace*{1pt}
the gradient of $\mathbf{f}$, that is, the matrix $(\nabla
\mathbf{f})_{i,j}:=\frac{\partial\mathbf{f}_i}{\partial x_j}$. We
will simply
write $(\nabla\mathbf{f}) \mathbf{g}$ for the column-vector $(\sum_j
\frac{\partial\mathbf{f}_i}{\partial x_j}\mathbf{g}_j )_i$ [instead
of the usual
``$(\mathbf{g}\cdot\nabla) \mathbf{f}$''].

\item[--] $C$ and $C(T)$ are finite positive constants that may
change from line to line.
\end{enumerate}

\section{The weak 3d-vortex equation and a probabilistic
interpretation of the external field}\label{sec2}

Let us recall a that vector field $w\dvtx\mathbb{R}^3\to\mathbb{R}^3$ with
components in $\mathcal{D}'$, and such that $\int_{\mathbb
{R}^3}\nabla
f(x)w(x)\,dx=0 $ for all $f\in\mathcal{D}$, is said to have \textit{null
divergence in the distribution sense}. We write it $\operatorname{div} w= 0$.

If the following two conditions hold, we shall say that $w_0\dvtx\mathbb
{R}^3\to
\mathbb{R}^3$ and
$\mathbf{g}\dvtx\mathbb{R}_+\times\mathbb{R}^3\to\mathbb{R}^3$
satisfy the hypothesis:

\begin{enumerate}
\item[$(\mathrm{H}_p)$:]\hypertarget{Hp}
\mbox{}
\begin{itemize}
\item there exists $p\in[1,\infty[$ such that
$w_0 \in(L^p(\mathbb{R}^3))^3$ and $ \mathbf{g}(t,\cdot) \in
(L^p(\mathbb{R}^3))^3$
for all $t\in[0,T]$, and
$ \sup_{t \in[0,T]} \|\mathbf{g}(t,\cdot
)\|_p<\infty$;
\item
$\operatorname{div} w_0=0$ and $\operatorname{div} \mathbf{g}(t,\cdot)=0$ for all $t\in[0,T]$.
\end{itemize}
\end{enumerate}

A necessary assumption for our probabilistic approach will be that
\hyperlink{Hp}{$(\mathrm{H}_p)$} holds with $p=1$. We then denote
\[
\|\mathbf{g}\|_{1,T}:=\int_0^T \int_{\mathbb{R}^3}|\mathbf{g}(s,x)|\,dx \,ds.
\]

In that functional setting, the following notion of solution to
(\ref{3DVF}) will appear to be natural:
\begin{definicion}\label{weakeq} Let $w_0$ and $\mathbf{g}$ satisfy
\hyperlink{Hp}{$(\mathrm{H}_1)$}. A function $\mathbf{w}\in
L^{\infty}([0,T],\break
(L^1(\mathbb{R}^3))^3)$ is a weak solution on $[0,T]$ of the vortex
equation with initial condition $w_0$ and external field $\mathbf{g}$ (or
``weak solution'') if:
\begin{longlist}
\item For $i,j,k=1,2,3$,
%
%
\begin{eqnarray}\label{integrability}
\int_{[0,T]\times\mathbb{R}^3} |\mathbf{w}_i(t,x)||\mathbf
{K}(\mathbf{w})_j(t,x)| \,dx \,dt &<&
\infty, \nonumber\\[-8pt]\\[-8pt]
\int_{[0,T]\times\mathbb{R}^3} |\mathbf{w}_i(t,x)| \biggl|\frac
{\partial\mathbf{K}(\mathbf{w}
)_j}{\partial x_k}(t,x) \biggr| \,dx\,
dt&<&\infty.\nonumber
\end{eqnarray}

\item For any $\mathbf{f}\in(C^{1,2}_b)^3$,
%
%
\begin{eqnarray}\label{weak}\quad
&&\int_{\mathbb{R}^3}\mathbf{f}(t,y)\mathbf{w}(t,y)\,dy \nonumber\\
&&\qquad=\int
_{\mathbb{R}^3} \mathbf
{f}(0,y)w_0(y)\,dy +
\int_0^t \int_{\mathbb{R}^3}\mathbf{f}(s,y)\mathbf{g}(s,y)\,dy\, ds
\nonumber\\[-8pt]\\[-8pt]
&&\qquad\quad{} + \int_0^t
\int_{\mathbb{R}^3} \biggl[ \frac{\partial\mathbf{f}}{\partial s}(s,y)
+\nu\triangle\mathbf{f}(s,y)\nonumber\\
&&\qquad\quad\hspace*{47.3pt}{}
+ \nabla\mathbf{f}(s,y) \mathbf{K}(\mathbf{w})(s,y)+\mathbf
{f}(s,y)\nabla
\mathbf{K}(\mathbf{w})(s,y) \biggr]\mathbf{w}(s,y) \,dy \,ds.\nonumber
\end{eqnarray}
\end{longlist}
\end{definicion}
\begin{rem} We observe that for any
function $\mathbf{v}\dvtx\mathbb{R}^3\to\mathbb{R}^3$ in $L^1$, the
functions $\mathbf{K}(\mathbf{v})$ and
$\nabla\mathbf{K}(\mathbf{v})$ are defined a.e. on $ \mathbb{R}^3$.
Indeed, the first
one can be bounded by a (scalar) Riesz potential operator (see
Stein \cite{Stein}), and thus belongs to a suitable weak Lebesgue
space. The second one is defined through a singular integral
operator acting on $\mathbf{v}$ (see, e.g., \cite{BerMa} for this
fact), and
this implies (see also \cite{Stein}) that it is an almost
everywhere defined function of some other weak Lebesgue space.
\end{rem}

We next introduce the central probabilistic objects we shall be
dealing with, which extend the ideas introduced in two dimensions
in \cite{FM}.
\begin{definicion}\label{margsP}
We write $\mathcal{C}_T:=[0,T]\times C([0,T],\mathbb{R}^3\times
\mathbb{R}^{3\otimes3})$. The canonical process in $\mathcal{C}_T$
will be
denoted by $(\tau, X,\Phi)$, and the space of probability measures
on $\mathcal{C}_T$ is written $\mathcal{P}(\mathcal{C}_T)$.

For an element $P\in\mathcal{P}(\mathcal{C}_T)$, we write
$P^{\circ}=\operatorname{law}(X)$ for the second marginal and $P'=\operatorname{law}(\Phi)$ for
the third marginal.
\end{definicion}

We shall also denote
%
%
\begin{eqnarray}\label{initlaw}
\bar{w}_0(x) & = & \frac{|w_0(x)|}{\|w_0\|_1+\|\mathbf{g}\|_{1,T}}\quad
\mbox{and } \nonumber\\[-8pt]\\[-8pt]
\bar{\mathbf{g}}(t,x)& = &\frac{|\mathbf{g}(t,x)|}{\|w_0\|_1+\|
\mathbf{g}\|_{1,T}}.\nonumber
\end{eqnarray}

We then define a probability measure $P_0(dt,dx)$ on $[0,T]\times
\mathbb{R}^3$ by
%
%
\begin{equation}\label{P0}
P_0(dt,dx) =\delta_0(dt) \bar{w}_0(x) \,dx + \bar{\mathbf{g}}(t,x) \,dx \,dt,
\end{equation}
together with the vectorial weight function
%
%
\begin{eqnarray}\label{h0}
h(t,x) &=& \mathbf{1}_{\{t=0\}}\frac{w_0(x)}{|w_0(x)|} (\|w_0\|_1+\|\mathbf
{g}\|
_{1,T} )\nonumber\\[-8pt]\\[-8pt]
&&{}
+\frac{\mathbf{g}(t,x)}{|\mathbf{g}(t,x)|} (\|w_0\|_1+\|\mathbf
{g}\|_{1,T} )\mathbf{1}_{\{t>0\}},\nonumber
\end{eqnarray}
where $\mathbf{1}$ denotes the indicator function and the convention
``$\frac{0}{0}=0$'' is made. We notice that
$|h(t,x)|=\|w_0\|_1+\|g\|_{1,T}$ or $0$. Moreover, we have
\begin{rem}\label{P0rem}
For measurable bounded functions $\mathbf{f}\dvtx[0,T]\times\mathbb
{R}^3\to\mathbb{R}^3$,
we have
\begin{eqnarray*}
&&\int_{[0,T]\times\mathbb{R}^3} \mathbf{f}(s,x)h(s,x)P_0(ds,dx)
\\
&&\qquad=
\int
_{\mathbb{R}^3}
\mathbf{f}(0,x)w_0(x)\,dx + \int_{[0,T]\times\mathbb{R}^3}
\mathbf
{f}(s,x)\mathbf{g}(s,x)\,dx
\,ds.
\end{eqnarray*}
\end{rem}

Consider now $Q\in\mathcal{P}(\mathcal{C}_T)$ such that for all
$\in[0,T]$, $\mathbb{E}^{Q}(|\Phi_t|)<\infty$.
Then, we can associate
with $Q$ a family of $\mathbb{R}^3$-valued vector measures
$(\tilde{Q}_t)_{t\in[0,T]}$ on $\mathbb{R}^3$, defined for all bounded
measurable function $\mathbf{f}\dvtx\mathbb{R}^3\to\mathbb{R}^3$ by
%
%
\begin{equation}\label{measvect}
\tilde{Q}_t(\mathbf{f})=\mathbb{E}^{Q}\bigl(\mathbf{f}(X_t)\Phi_t
h(\tau
,X_0)\mathbf{1}_{\{\tau\leq
t\}}\bigr).
\end{equation}

Moreover, $\tilde{Q}_t$ is absolutely continuous with respect to
$Q^{\circ}_t$, with
%
%
\begin{equation}\label{density}
\frac{d \tilde{Q}_t}{d Q^{\circ}_t}(x)=E^Q\bigl(\Phi_t
h(\tau,X_0)\mathbf{1}_{\{\tau\leq t\}}\vert X_t=x\bigr),
\end{equation}
and its total mass is bounded by $ (\|w_0\|_1
+\|\mathbf{g}\|_{1,T})\mathbb{E}^{Q}(|\Phi_t|) $.
\begin{definicion} We denote by $\mathcal{P}_b(\mathcal{C}_T)$ the
subset of
probability measures $Q\in\mathcal{P}(\mathcal{C}_T)$ under which
the process $\Phi$ belongs to $L^{\infty}( [0,T]\times\Omega
,dt\otimes Q)$.
\end{definicion}

Then, we consider the following nonlinear martingale problem:

\begin{enumerate}[]
\item[(MP):]\hypertarget{MP} to find $P\in\mathcal{P}_b(\mathcal{C}_T)$
such that:
\begin{itemize}
\item
$X_t=X_0$ in $[0,\tau]$, $P$-almost surely.

\item The law of $(\tau,X_0)$ under $P$
is $P_0$ given by (\ref{P0}), and $\tilde{P}_t$
constructed
according to (\ref{measvect}) has a bi-measurable
density family
$(t,x)\mapsto\tilde{ \rho}(t,x)$.

\item $f(t,X_t)-f (0,X_0)-\int_{0}^t \frac{\partial
f}{\partial s}(s,X_s)+ [\nu\triangle f (s,X_s)
+ \mathbf{K}(\tilde{\rho})(s,X_s)\nabla f(s,\break X_s) ] \mathbf
{1}_{ s\geq\tau
}\,ds$,
$0\leq t\leq T$, is a continuous $P$-martingale for all $f
\in\mathcal{C}_b^{1,2}$
w.r.t. the filtration
$\mathcal{F}_t=\sigma(\tau,(X_s,\Phi_s),s\leq t)$.

\item $\Phi_t=I_3 +\int_{0}^t \nabla
\mathbf{K}(\tilde{\rho})(s,X_s)\Phi_s \mathbf{1}_{ s\geq\tau}
\,ds$, for all
$0\leq t\leq T$, $P$-almost surely.
\end{itemize}
\end{enumerate}

The following statement partially explains the relation between
\hyperlink{MP}{($\mathrm{MP}$)} and (\ref{3DVF}), and will be useful later on:
\begin{lema}\label{MPtoweakeq}
Assume that the problem \textup{\hyperlink{MP}{($\mathrm{MP}$)}} has a solution $P \in\mathcal
{P}_b(\mathcal{C}_T)$ satisfying
%
%
\begin{equation}\label{integKrho}
E \biggl(\int_0^T | \mathbf{K}(\tilde{\rho})(t,X_t)|\,dt \biggr)
<\infty
\end{equation}
and
%
%
\begin{equation}\label{integgradKrho}
E \biggl(\int_0^T |\nabla
\mathbf{K}(\tilde{\rho})(t,X_t)|\,dt \biggr)<\infty.
\end{equation}
Then, $\tilde{\rho}$ is
a weak solution of the vortex equation with external force field
(\ref{weak}).
\end{lema}
\begin{pf}
The assumptions on $P$ imply that point (i) in Definition
\ref{weakeq} is satisfied and, moreover, that $\int_0^t
\mathbf{K}(\tilde{\rho})(s,X_s)\,ds$ and $\int_0^t \nabla
\mathbf{K}(\tilde{\rho})(s,X_s)\,ds$ are both processes with integrable
variation (and thus absolutely continuous on $[0,T]$). Since under
$P$ the process $\Phi_t$ is almost surely bounded in $[0,T]$, it
follows that it has finite variation too.

On the other hand, the martingale associated with $f \in\mathcal
{C}_b^{1,2}$ in \hyperlink{MP}{($\mathrm{MP}$)} equals
\begin{eqnarray*}
&&f(t,X_t)-f (\tau\wedge t,X_0)\\
&&\qquad{} -\int_{0}^t \biggl[\frac{\partial
f}{\partial s}(s,X_s)+\nu\triangle f (s,X_s) +
\mathbf{K}(\tilde{\rho})(s,X_s)\nabla f(s,X_s) \biggr] \mathbf{1}_{
s\geq\tau}\,ds
\end{eqnarray*}
thanks to the first condition of \hyperlink{MP}{($\mathrm{MP}$)}.

Therefore, by It\^{o}'s product rule, we see that for each $\mathbf
{f}\in
(C^{1,2}_{b})^3$
\begin{eqnarray*}
&&\mathbf{f}(t,X_t)\Phi_t-\mathbf{f}(\tau\wedge,X_0)\\
&&\qquad{}-\int_0^t
\biggl[\frac{\partial
\mathbf{f}}{\partial s}(s,X_s)+\nu\triangle\mathbf{f}(s,X_s) +
\nabla\mathbf{f}(s,X_s)\mathbf{K}(\tilde{\rho})(s,X_s)\\
&&\qquad\hspace*{138.4pt}{} + \mathbf{f}(s,X_s)
\nabla
\mathbf{K}(\tilde{\rho})(s,X_s) \biggr]\Phi_s \mathbf{1}_{\{s\geq
\tau\}}\,ds
\end{eqnarray*}
is a local martingale issued from $0$. Moreover, the assumptions
(\ref{integgradKrho}) and (\ref{integKrho}) on $\tilde{\rho}$ and
the fact that $\Phi$ is bounded imply that it is a true
martingale. Consequently, as $h(\tau,X_0)\mathbf{1}_{\{\tau\leq t\}
}$ is
$\mathcal{F}_0$-measurable and $\mathbf{1}_{\{\tau\leq s\}\cap\{
\tau\leq
t\}}=\mathbf{1}_{\{\tau\leq s\}}$ for $s\leq t$, we see that
%
%
\begin{eqnarray}
&&
E^P  \bigl( \mathbf{f}(t,X_t) \Phi_t h(\tau,X_0)\mathbf{1}_{\{
\tau\leq
t\}} \bigr) -E^P \bigl( \mathbf{f}(\tau,X_0)h(\tau,X_0)\mathbf
{1}_{\{\tau
\leq t\}} \bigr) \nonumber\\
&&\qquad{} - E^P \biggl( \int_0^t \biggl[\frac{\partial\mathbf{f}}{\partial
s}(s,X_s)+\nu\triangle\mathbf{f}(s,X_s)\nonumber\\[-8pt]\\[-8pt]
&&\qquad\hspace*{52.4pt}{} +\nabla
\mathbf{f}(s,X_s)\mathbf{K}(\tilde{\rho})(s,X_s) \nonumber\\
&&\qquad\hspace*{52.4pt}{}
+ \mathbf{f}(s,X_s)
\nabla\mathbf{K}(\tilde{\rho})(s,X_s) \biggr]\Phi_s
h(\tau,X_0)\mathbf{1}_{\{\tau\leq
s\}}\,ds \biggr)=0.\nonumber
\end{eqnarray}
Recalling that $\tilde{\rho}(t)$ is the
density of the vector measure (\ref{measvect}) for $Q=P$, the first
term in the previous equation is seen
to be equal to $\int\mathbf{f}(t,x)\tilde{\rho}(t,x)\,dx$. The second
term is equal to the expression in Remark \ref{P0rem} with
$\mathbf{f}(s,x)$ replaced by $\mathbf{f}(s,x)\mathbf{1}_{s\leq t}$,
that is,
$\int\mathbf{f}(0,y)w_0(y)\,dy +
\int_0^t \int\mathbf{f}(s,y)\mathbf{g}(s,y)\,dy \,ds$. The third expectation
can be
interchanged with the time integral
thanks
to the assumptions and Fubini's theorem, and the result follows
using again the definition of $\tilde{\rho}(s)$ in the resulting
time integral.
\end{pf}

The proof of the well posedness of problem \hyperlink{MP}{($\mathrm{MP}$)} will be based
on analytical results about the ``mild form'' of the vortex
equation (\ref{3DVF}), which we state in next section. These will
in particular provide a framework where the conditions required in
Lemma \ref{MPtoweakeq} will hold.

\section{The mild vortex equation in $L^p$ with an external field}
\label{sec3}

We shall next introduce the mild formulation of the forced vortex
equation. We refer the reader to the book of Lemari\'{e}-Rieusset
\cite{LMR} for a comprehensive account on the mild-form approach
to the Navier--Stokes equation in its velocity form. Our techniques
are adapted from that framework.

We denote the heat kernel in $\mathbb{R}^3$ by
%
%
\begin{equation}\label{heatkernel}
G^{\nu}_t (x):=(4\pi\nu t)^{-{3/2}} \exp\biggl(-
\frac{|x|^2}{4\nu t} \biggr),
\end{equation}
where $\nu>0$. One has
\begin{lema}\label{estimsw}
For all $p\in[1,\infty]$, $r\geq p$ and $w\in(L^p)^3$, there
exist positive constants $\bar{C}_0(p;r)$ and
$\bar{C}_1(p;r)$ such that for all $t>0$:
\begin{longlist}
\item$\|G_t^{\nu}*w\|_r\leq\bar{C}_0(p;r)t^{-{3/2}
({1/p}-{1/r})} \|w\|_p$,
\item$ \|\nabla
G_t^{\nu}*w\|_r\leq\bar{C}_1(p;r)t^{-{1/2}-{3/2}
({1/p}-{1/r})} \|w\|_p$.
\end{longlist}
\end{lema}
\begin{pf} Use Young's inequality and the well-known estimates
\[
\sup_{t\geq0} \|G^{\nu}_ t\|_m t^{{3/2} -
{3/(2m)}}<\infty, \qquad\sup_{t\geq0} \|\nabla G^{\nu}_
t\|_m t^{2-{3/(2m)}} <\infty.
\]
\upqed\end{pf}
\begin{definicion} Let $w_0$ and $\mathbf{g}$ be functions satisfying
\hyperlink{Hp}{$(\mathrm{H}_p)$} for some $p \in[1,\infty]$. A function $\mathbf
{w}\in
L^{\infty}([0,T], (L^p(\mathbb{R}^3))^3)$ is a mild solution on $[0,T]$
of the vortex equation with initial condition $w_0$ and external
field (or ``mild solution'') if:
\begin{longlist}
\item The functions
$\mathbf{K}(\mathbf{w})_i(t,x):=\mathbf{K}(\mathbf{w}(t,\cdot
))_i(x)$, $i=1,2,3$ are defined a.e.
on $[0,T]\times\mathbb{R}^3$ and satisfy the integrability conditions
(\ref{integrability}).
\item For $dt$-almost every $t$, the following identity
holds in $(L^p)^3$:
%
%
\begin{eqnarray}\label{eqmild}\qquad
\mathbf{w}(t,x)&=&G^{\nu}_t* w_0(x)+ \int_0^t G^{\nu}_{t-s}*\mathbf
{g}(s,\cdot)
(x)\,ds \nonumber\\
&&{}+ \sum_{j=1}^3\int_0^t \int_{\mathbb{R}^3}\frac{\partial
G^{\nu}_{t-s}}
{\partial y_j}
(x-y) [\mathbf{K}(\mathbf{w})_j(s,y)\mathbf{w}(s,y)
\\
&&\hspace*{127.4pt}{}- \mathbf{w}_j(s,y)\mathbf{K}(\mathbf{w})(s,y) ]\,dy \,ds.\nonumber
\end{eqnarray}
\end{longlist}
\end{definicion}

We shall state in Theorems \ref{exist1} and \ref{regularity} below
the analytical results we need about (\ref{eqmild}). As we shall
see, that equation will admit an abstract formulation which is the
same as in the case $\mathbf{g}=0$, and so we will be able to adapt the
techniques in \cite{F1} with no difficulties. We therefore
provide an abbreviated account of these results.

We shall simultaneously deal with a family of ``mollified''
versions of (\ref{eqmild}). Consider a smooth function
$\varphi\dvtx\mathbb{R}^3\to\mathbb{R}$ satisfying:
\begin{longlist}
\item$\int_{\mathbb{R}^3}\varphi(x)\,dx=1$,
\item${\int_{\mathbb{R}^3}}|x| |\varphi(x)|\,dx<\infty$,
\end{longlist}
which is called a ``cutoff function of order $1$.'' For
$\varepsilon>0$, let $\varphi_{\varepsilon}\dvtx\mathbb{R}^3\to\mathbb
{R}$ denote
the regular approximation of the Dirac mass
$\varphi_{\varepsilon}(x)=\frac{1}{\varepsilon^3}\varphi(\frac
{\varepsilon}{x})$.
We define the convolution operators
%
%
\begin{equation}\label{biotsavop}
\mathbf{K}^{\varepsilon}(w)(x):=\int_{\mathbb{R}^3}K_{\varepsilon}
(x-y)\wedge w(y) \,dy,
\end{equation}
where $K_{\varepsilon}:=\varphi_{\varepsilon}*K
=\mathbf{K}(\varphi_{\varepsilon})$. The fact that $K_{\varepsilon
}$ is a
regular function will follow from part (ii) in Lemma
\ref{contK} below. To unify notation, we also write $K_0=K$ and $
\mathbf{K}^{0}(w)(x):=\mathbf{K}(w)(x)$.

We introduce the
family $\{\mathbf{B}^{\varepsilon}\}_{\varepsilon\geq0} $ of operators
(formally) defined on functions $\mathbf{w},\mathbf{v}\dvtx[0,T]\times
\mathbb{R}^3 \to
\mathbb{R}^3$ by
%
%
\begin{eqnarray}\label{bil}\quad
&&\mathbf{B}^{\varepsilon}(\mathbf{w},\mathbf{v})(t,x)
\nonumber\\
&&\qquad=\int_0^t
\sum_{j=1}^3\int_{\mathbb{R}
^3}\frac{\partial G^{\nu}_{t-s}}{\partial y_j}
(x-y) \\
&&\qquad\quad\hspace*{46.2pt}{}\times [\mathbf{K}^{\varepsilon}(\mathbf{w})_j(s,y)\mathbf{v}(s,y)-
\mathbf{v}_j(s,y)\mathbf{K}^{\varepsilon}(\mathbf{w})(s,y) ]\,dy
\,ds.\nonumber
\end{eqnarray}

We are interested in the following family of ``abstract''
equations, for $\varepsilon\geq0$:
%
%
\begin{equation}\label{abs}
\mathbf{v}=\mathbf{w}_0+ \mathbf{B}^{\varepsilon}(\mathbf
{v},\mathbf{v}),
\end{equation}
where
\[
\mathbf{w}_0(t,x):=G^{\nu}_t* w_0(x)+ \int_0^t G^{\nu
}_{t-s}*\mathbf{g}(s,\cdot)
(x) \,ds.
\]

For a given time interval $[0,T]$ we shall work in the Banach
spaces
\[
\mathbf{F}_{0,r,(T;p)},\qquad \mathbf{F}_{1,r,(T;p)},\qquad \mathbf
{F}_{0,p,T}\quad\mbox{and}\quad \mathbf{F}_{1,p,T}
\]
with norms, respectively, defined by:
\begin{itemize}
\item$ {\tn\mathbf{w}\tn_{0,r,(T;p)}:=\sup_{0\leq t
\leq T}
t^{{3/2}({1/p}-{1/r})}\|\mathbf{w}(t)\|_r }$,
\item
$ \tn\mathbf{w}\tn_{1,r,(T;p)}:=\sup_{0\leq t \leq T}
\{ t^{{3/2}({1/p}-{1/r})}\|\mathbf{w}(t)\|
_r +
t^{{1/2}+{3/2}({1/p}-{1/r})}\times\break
{\sum_{k=1}^3} \| \frac{\partial\mathbf{w}(t)}{\partial
x_k} \|_r
\} $,
\item$ {\tn\mathbf{w}\tn_{0,p,T}:=\tn\mathbf{w}
\tn_{0,p,(T;p)}}$ and
\item$ {\tn\mathbf{w}
\tn_{1,p,T}:=\tn\mathbf{w}\tn_{1,p,(T;p)}}$.
\end{itemize}

The following continuity property of the Biot--Savart kernel is
crucial:
\begin{lema}\label{contK}
Let $1< p <3$ be given and $q\in(\frac{3}{2},\infty)$ be defined
by $\frac{1}{q}=\break\frac{1}{p}-\frac{1}{3}$.
\begin{longlist}
\item For every $w \in(L^3)^p$, the integral (\ref{biotsavop})
is absolutely convergent for almost every $x$ and one has
$\mathbf{K}^{\varepsilon}(w)\in(L^q)^3$. There exists further a
positive constant $\tilde{C}_{p,q}$ such that
%
%
\begin{equation}\label{Kpq}
\sup_{\varepsilon\geq0}\| \mathbf{K}^{\varepsilon}(w)\|_q\leq
\tilde{C}_{p,q} \|w\|_p
\end{equation}
for all $w \in(L^p)^3$.
\item If moreover $w \in
(W^{1,p})^3$, then we have
$\mathbf{K}^{\varepsilon}(w)\in(W^{1,q})^3$, with\break $\frac{\partial}
{\partial x_k}\mathbf{K}^{\varepsilon}(w)=\mathbf{K}^{\varepsilon}
(\frac{\partial w }{\partial x_k} )$, and
%
%
\begin{equation}\label{KWpq}
\sup_{\varepsilon\geq0} \biggl\| \frac{\partial
\mathbf{K}^{\varepsilon}(w)}{\partial x_k} \biggr\|_q\leq
\tilde{C}_{p,q} \biggl\|\frac{\partial w }{\partial
x_k} \biggr\|_p
\end{equation}
for all $k=1,2,3$.
\end{longlist}
\end{lema}
\begin{pf} See
Lemma 2.2 in \cite{F1} for the case $\varepsilon=0$ and Remark 4.3
therein for the general case.
\end{pf}
\begin{lema}\label{contB}
\textup{(i)}
Let $p\in[1,3)$ and assume \hyperlink{Hp}{$(\mathrm{H}_p)$}. Then, we
have for all $r\in[p,\frac{3p}{3-p})$ that
\[
\mathbf{w}_0\in F_{1,r,(T;p)}\qquad \mbox{with } \tn\mathbf{w}_0
\tn_{1,r,(T;p)}\leq C(r,p)(\|w_0\|_p + T \tn\mathbf{g}\tn_{0,p,T})
\]
for some finite constant $C(r,p)>0$.

{\smallskipamount=0pt
\begin{longlist}[(ii)]
\item[(ii)] Let $\frac{3}{2}< p<3, p\leq l <\min\{\frac{6p}{6-p},3\}$
and $ \frac{3l}{6-l}\leq
l'< \frac{3l}{6-2l}$. Then, there exists a finite constant
$C_1(l,l';p)$ not depending on $T>0$ such that for all $\mathbf
{w},\mathbf{v}\in
\mathbf{F}_{1,l,(T;p)}$,
\[
\sup_{\varepsilon\geq0}\tn\mathbf{B}^{\varepsilon}(\mathbf
{w},\mathbf{v}) \tn
_{1,l',(T;p)} \leq C_1(l,l';p)
T^{1-{3}/({2p})} \tn\mathbf{w}
\tn_{1,l,(T;p)} \tn\mathbf{v}\tn_{1,l,(T;p)},
\]
where
$1-\frac{3}{2p}>0$.
\end{longlist}}
\end{lema}
\begin{pf}
Part (i) follows from Lemma \ref{estimsw}. To
bound the time integral we use, moreover, the fact that for all
$r\geq p$, on has
\[
\biggl\|\int_0^t
G^{\nu}_{t-s}*\mathbf{g}(s,\cdot) \,ds \biggr\|_r\leq C
t^{1+{1/r}-{1/p}} \Bigl({\sup_{t\in
[0,T]}}\|g_t\|_p \Bigr).
\]
On the other hand, since $t\mapsto
t^{-{1/2}+{3/2}({1/r}-{1/p})}$ is
integrable in $0$ if and only if $r<\frac{3p}{3-p}$, we have
\[
\biggl\|\nabla\biggl(\int_0^t
G^{\nu}_{t-s}*\mathbf{g}(s,\cdot) \,ds \biggr) \biggr\|_r\leq C'
t^{{1/2}+{3/2}({1/r}-{1/p})} \Bigl({\sup
_{t\in
[0,T]}}\|\mathbf{g}(t,\cdot) \|_p \Bigr)
\]
from where the statement follows. Part (ii) uses Lemma
\ref{contK} and is proved in parts (ii) and (iv) of
Proposition 3.1 in \cite{F1}. See also Remarks 4.3 and 6.3 therein
for the uniformity (in $\varepsilon\geq0 $) of the bounds.
\end{pf}
\begin{rem}
Observe that the previous lemma, in particular, implies
(taking $p=r=l=l'$) that for $p\in(\frac{3}{2},3)$, the abstract
equation (\ref{abs}) makes sense in $\mathbf{F}_{1,p,T}$ for each
$\varepsilon\geq0$.
\end{rem}

Now we can state the extension of Theorem 3.1 in \cite{F1} to the
3d-vortex equation with external field.
\begin{teorema}\label{exist1}
Assume that \hyperlink{Hp}{$(\mathrm{H}_p)$} for some $\frac{3}{2}<p<3$.
\begin{enumerate}[(a)]
\item[(a)] For each $T>0$ and $\varepsilon\geq0$, equation (\ref
{abs}) has, at most,
one solution in $\mathbf{F}_{0,p,T}$.

\item[(b)] There is a constant
$\Gamma_0(p)>0$ independent of $\varepsilon\geq0$ such that for
all $T>0$, $w_0 $ and $\mathbf{g}$
satisfying
\[
T^{1-{3/(2p)}} (\|w_0\|_p+T\tn\mathbf{g}
\tn_{0,p,\theta} )<\Gamma_0(p),
\]
each one of (\ref{abs}) with $\varepsilon\geq0$, has a
solution $\mathbf{w}^{\varepsilon}\in\mathbf{F}_{1,p,T}$. Moreover,
we have
\[
\sup_{\varepsilon\geq0}\tn\mathbf{w}^{\varepsilon}\tn
_{1,p,T}\leq2 \tn
\mathbf{w}_0 \tn_{0,p,T}.
\]
\end{enumerate}
\end{teorema}
\begin{pf} For later purposes, we give, in detail, the argument of
\cite{F1}. By Lem\-ma~\ref{estimsw}(ii) (with $p$ in the place of $r$ and
$\frac{3p}{6-p}$ in that of $p$) and Lemma \ref{contK}(i),
we have for all $\mathbf{v},\mathbf{w}\in\mathbf{F}_{0,p,T}$ that
\[
\| \mathbf{B}^{\varepsilon}(\mathbf{w},\mathbf{v})(t)\|_p\leq C\int_0^t
(t-s)^{-{3}/({2p})} \|\mathbf{w}(s)\|_p\|\mathbf{v}(s)\|_p \,ds.
\]
It follows that
if $\mathbf{w}$ and $\mathbf{v}$ are two solutions, one has
\[
\|\mathbf{w}(t) -\mathbf{v}(t)\|_p \leq C ( \tn\mathbf{w}\tn
_{0,p,T}+\tn\mathbf{v}
\tn_{0,p,T} )\int_0^t
(t-s)^{-{3}/({2p})} \|\mathbf{w}(s)-\mathbf{v}(s)\|_p
\,ds
\]
and iterating the latter sufficiently many times [using the
identity $ \int_0^t s^{\varepsilon-1}(t-s)^{\theta-1} \,ds =C
t^{\varepsilon+ \theta-1}$ for $\theta,\varepsilon>0$] we get
$\|\mathbf{w}(t) -\mathbf{v}(t)\|_p \leq C \int_0^t \|\mathbf
{w}(s)-\mathbf{v}(s)\|_p \,ds$.
Gronwall's lemma concludes the proof.

(b) We
notice that for $T>0$ small enough, one has
\[
4 C(p,p) C_1(p,p;p)
T^{1-{3}/({2p})} (\|w_0\|_p + T \tn\mathbf{g}\tn_{0,p,T})<1,
\]
where $C(p,p)$ and $C_1(p,p;p)$ are, respectively, the constants in
parts (i) and (ii) of Lemma \ref{contB} with all
parameters equal to $p$. From this and Lemma \ref{contB}(i),
the same contraction argument used in Theorem 3.1(b) of
\cite{F1} can be applied here in the space $\mathbf{F}_{1,p,T}$.
\end{pf}

We observe that for $\mathbf{v}\in\mathbf{F}_{0,p,T}$, with $p\in
(\frac{3}{2},3)$ we have $\mathbf{K}(\mathbf{v}) \in\mathbf
{F}_{0,q,T}$ for $q\in
(3,\infty)$. The previous global uniqueness and local existence
result also holds in that space, and one can, moreover, show that
the solution $\mathbf{w}(t)\in(L^p)$ is a continuous function of $t$.
That type of result corresponds to a ``vorticity version'' of
Kato's theorem for the mild Navier--Stokes equation in $(L^q)^3$,
$q\in(3,\infty)$ (see~\cite{LMR}, Theorem 15.3(A)).

We shall, later on, need additional regularity properties of the
function $\mathbf{w}^{\varepsilon}$ and, more importantly, their uniformity
in $\varepsilon\geq0$. These results will rely on continuity
properties of the ``derivative'' of the Biot--Savart operator.
\begin{lema}\label{contgradK}
Let $1<r<\infty$.
\begin{longlist}
\item For all\vspace*{1pt} $w \in(L^r)^3$ and $\varepsilon\geq0$, we have
$\frac{\partial
}{\partial x_k}\mathbf{K}^{\varepsilon}(w)\in(L^r)^3$ for
$k=1,2,3$. There exists further a positive constant $C_r$
depending only on $r$ such that
%
%
\begin{equation}\label{HLp}
\sup_{\varepsilon\geq0} \biggl\| \frac{\partial
\mathbf{K}^{\varepsilon}(w)_j}{\partial x_k} \biggr\|_r\leq
\tilde{C}_r \|w\|_r
\end{equation}
for all $j=1,2,3$, where $\mathbf{K}^{\varepsilon}(w)_j$ is the
$j$th component of $\mathbf{K}^{\varepsilon}(w)$.
\item If,
moreover, $w \in(W^{1,r})^3$, we then have
$\frac{\partial}{\partial
x_k}\mathbf{K}^{\varepsilon}(w)\in(W^{1,r})^3$, with
$\frac{\partial}{\partial x_i} (\frac{\partial}{\partial
x_k}\mathbf{K}^{\varepsilon}(w) )= \frac{\partial}{\partial
x_k}\mathbf{K}^{\varepsilon}(\frac{\partial}{\partial x_i} w) $
and
%
%
\begin{equation}\label{HW1p}
\sup_{\varepsilon\geq0} \biggl\| \frac{\partial^2
\mathbf{K}^{\varepsilon}(w)_j}{\partial x_i\,\partial x_k} \biggr\|_r
\leq\tilde{C}_r \biggl\|\frac{\partial w }{\partial
x_i} \biggr\|_r
\end{equation}
for all $i,k=1,2,3$.
\end{longlist}
\end{lema}
\begin{pf} See Lemma 3.1 and Remark 4.3 in \cite{F1} for the proof,
which relies on the fact that $\mathbf{w}\mapsto\frac{\partial
\mathbf{K}(w)}{\partial x_k}$ is a singular integral operator.
\end{pf}
\begin{teorema}\label{regularity}
For $p\in(\frac{3}{2},3)$, let $\mathbf{w}^{\varepsilon}\in\mathbf
{F}_{1,p,T}$,
$\varepsilon\geq0$ be the solution of (\ref{abs}) given by
Theorem \ref{exist1}, and write
$\mathbf{u}^{\varepsilon}(s,x):=\mathbf{K}^{\varepsilon}(\mathbf
{w}^{\varepsilon})(s,x)$.
Let $\mathcal{C}^{\alpha}$ denote the space of H\"{o}lder continuous
functions $\mathbb{R}^3\to\mathbb{R}^3$ of index $\alpha\in(0,1)$.
\begin{longlist}
\item For all $r\in[p,\frac{3p}{3-p})$, we have
\[
\sup_{\varepsilon\geq0} \tn\mathbf{w}^{\varepsilon} \tn_{1,r,(T;p)}
<\infty.
\]
\item
We have
%
%
\begin{equation}\label{estimu}
\sup_{\varepsilon\geq0}\sup_{t\in[0,T]} t^{{1/2}}
\{\|\mathbf{u}^{\varepsilon}(t)\|_{\infty}+\|\mathbf
{u}^{\varepsilon}(t)\|
_{\mathcal{C}^{({2p-3})/{p}}} \} <\infty.
\end{equation}
\item For all $r\in(3,\frac{3p}{3-p})$, $i=1,2,3$ we have
%
%
\begin{equation}\label{estimgradu}
\sup_{\varepsilon\geq0}\sup_{t\in[0,T]}
t^{{1/2}+{3/2}({1/p}-{1/r})} \biggl\{
\biggl\|\frac{\partial\mathbf{u}^{\varepsilon}(t)}{\partial
x_i} \biggr\|_{\infty}+ \biggl\|\frac{\partial\mathbf
{u}^{\varepsilon
}(t)}{\partial
x_i} \biggr\|_{\mathcal{C}^{1-{3/r}}} \biggr\} <\infty.
\end{equation}
In particular, the functions
\[
t\mapsto\|\mathbf{u}(t)\|_{\infty} \quad\mbox{and}\quad
t\mapsto\biggl\|\frac{\partial\mathbf{u}(t)}{\partial
x_i} \biggr\|_{\infty},\qquad i=1,2,3,
\]
belong to $L^1([0,T],\mathbb{R})$.
\end{longlist}
\end{teorema}
\begin{pf} Observe that
parts (i) and (ii) of Lemma \ref{contB} provide an estimate
of the form
\[
\tn
\mathbf{w}^{\varepsilon}\tn_{1,l',(T;p)}\leq C(l',p)(\|w_0\|_p+ T\tn
\mathbf{g}
\tn_{0,p,T}) +\Lambda(T,l,l')A_l^2
\]
for suitable $l$ and $l'$
and with $\Lambda(T,l,l')$ a uniform upper bound for the norms of
the operators $\mathbf{B}^{\varepsilon}\dvtx(\mathbf
{F}_{1,l,(T;p)})^2\to
\mathbf{F}_{1,l',(T;p)}$ and $A_l$ a given upper bound
of $\tn\mathbf{w}^{\varepsilon}\tn_{1,l,(T;p)}$. Then, starting
from the fact
that the functions
$\mathbf{w}^{\varepsilon}\in\mathbf{F}_{1,p,(T;p)}=\mathbf
{F}_{1,p,T}$ are uniformly
bounded in $\varepsilon\geq0$, we can apply several times Lemma
\ref{contB} and the previous
inequality (using, also, the fact that $\mathbf{w}_0 \in\mathbf
{F}_{1,l',(T;p)}$
for all $l' \in[p,\frac{3p}{3-p})$), and obtain an increasing
sequence $l'=l_n$ such that $l_0=p$, $l_n \nearrow
\frac{3p}{3-p}$, and
$\mathbf{w}^{\varepsilon}\in\mathbf{F}_{1,l_n,(T;p)}$ with $\tn
\mathbf{w}^{\varepsilon
}\tn
_{1,l_n,(T;p)}$
controlled in terms of $\tn\mathbf{w}^{\varepsilon}\tn_{1,l_{n-1},(T;p)}$
and $\tn\mathbf{w}_0\tn_{1,l_n,(T;p)}$. One can thus chose $N$ large enough
such that
$l_N\geq r$ and conclude with an
interpolation inequality in the spaces $\mathbf{F}_{1,l,(T;p)}$.
We refer to the proof of Theorem 3.2(ii) in \cite{F1} for this and
for an explicit construction of the sequence~$l_n$.

Next, Lemma \ref{contK} and Theorem \ref{exist1} imply that for
$q=\frac{3p}{3-p}>3$,
\[
\sup_{\varepsilon\geq0}\tn\mathbf{u}^{\varepsilon}
\tn_{{1,q,T}}\leq C\sup_{\varepsilon\geq0}\tn\mathbf
{w}^{\varepsilon}
\tn_{{1,p,T}} \leq C '(\|w_0\|_p+ T\tn\mathbf{g}
\tn_{0,p,T}) .
\]
Using the continuous embedding of $(W^{1,m})^3$
into $(L^{\infty})^3\cap\mathcal{C}^{1-{3}/{m}}$ for all $m>3$,
we deduce part (ii), taking $m=q$. To prove part (iii)
we use part (i), Lemma \ref{contgradK} and the same embedding
result as before but with $m=r$. See Corollary 3.1 in \cite{F1}
for details.
\end{pf}

\section{The nonlinear process}\label{sec4}

We shall, in this section, use the notation $F_{0,p,T}$, $F_{1,p,T}$,
$F_{0,r,(T;p)}$ and $F_{1,r,(T;p)}$ for the scalar-function
analogues of the spaces $\mathbf{F}$ defined in Section
\ref{sec3}.\vadjust{\goodbreak}

We also need the following definition.
\begin{definicion} $\mathcal{P}_{b,{3}/{2}}^T$ is the space of
probability measures $Q\in
\mathcal{P}_b(\mathcal{C}_T)$ satisfying the following conditions:
\begin{itemize}
\item For each $t\in[0,T]$, $Q^{\circ}_t(dx)$ defined in
Definition \ref{margsP} is
absolutely continuous with respect to Lebesgue's measure.

\item The family of densities of $(Q^{\circ}_t(dx))_{t\in[0,T]}$,
which we denote by $(t,x)\mapsto
\rho^Q(t,x)$, has a version that
belongs to $F_{0,p,T}$ for some $p>\frac{3}{2}$.

\item The family of densities of the vectorial measures $(\tilde
{Q}_t(dx))_{t\in[0,t]}$ [cf. (\ref{measvect})], which
we denote by $(t,x)\mapsto\tilde{\rho}^Q(t,x)$, satisfies
$\operatorname{div}
\tilde{\rho}_t^Q=0$ for $dt$-almost every $t\in[0,T]$.
\end{itemize}
\end{definicion}

We are ready to study the nonlinear process described in
\hyperlink{MP}{($\mathrm{MP}$)}.
\begin{teorema}\label{teoMP} Assume that \hyperlink{Hp}{$(\mathrm{H}_1)$} and
\hyperlink{Hp}{$(\mathrm{H}_p)$} are satisfied for some
$p\in(\frac{3}{2},3)$. Then, the following hold:
\begin{enumerate}[(a)]
\item[(a)] For every $T>0$, the nonlinear martingale problem \hyperlink{MP}{($\mathrm{MP}$)}
has, at most, one solution $P$ in the class $\mathcal{P}_{b,
{3/2}}^T$. Moreover, if such a solution $P$ exists,
then the function defined by
\[
\mathbf{w}(t,x):=\tilde{\rho}{}^P(t,x)=\rho^P (t,x)E^P\bigl(\Phi_t
h(\tau,X_0)\mathbf{1}_{\{t\geq\tau\}} \vert X_t=x\bigr)
\]
is the unique solution in $\mathbf{F}_{0,1,T}\cap\mathbf{F}_{0,p,T}$
of the
mild equation (\ref{eqmild}).
\item[(b)] In a given filtered probability space
$(\Omega,\mathcal{F},\mathcal{F}_t,\mathbb{P})$, consider a standard
three-dimensional Brownian motion $B$, and an $\mathcal{F}_0$-measurable
random variable $(\tau,X_0)$ independent of $B$
with law $P_0$ [defined as in (\ref{P0})]. Then, on each interval
$[0,T]$, the McKean nonlinear stochastic differential equation
%
%
\begin{eqnarray}\label{nonlinSDE}
\mbox{\textup{(i)} }&& X_t=X_0+\sqrt{2\nu} \int_0^t \mathbf{1}_{\{s\geq\tau\}
} \,dB_s
+\int_0^t \mathbf{K}(\tilde{\rho})
(s,X_s)\mathbf{1}_{\{s\geq\tau\}}\,ds ,\nonumber\\
\mbox{\textup{(ii)} }&& \Phi_t=I_3+\int_0^t \nabla\mathbf{K}(\tilde{\rho
})(s,X_s)\Phi_s \mathbf{1}_{\{s\geq\tau\}} \,ds , \\
\mbox{\textup{(iii)} }&& \operatorname{law}(\tau,X,\Phi) \in\mathcal{P}_{b,{3}/{2}}^T \quad
\mbox{and}\quad \tilde{\rho}(t,x)=\tilde{\rho}^{\operatorname{law}(\tau,X,\Phi)}(t,x),\nonumber
\end{eqnarray}
has, at most, one pathwise solution. Moreover, if a solution exists,
its law is a solution of \hyperlink{MP}{($\mathrm{MP}$)}. Thus, by \textup{(a)},
uniqueness in law for (\ref{nonlinSDE}) holds.
\item[(c)] If the condition
\[
T^{1-{3}/({2p})} (\|w_0\|_p+T\tn\mathbf{g}
\tn_{0,p,\theta} )<\Gamma_0(p)
\]
is satisfied,\vspace*{1pt} where
$\Gamma_0(p)>0$ is the constant provided by Theorem \ref{exist1},
then a unique solution $P\in\mathcal{P}_{b,{3/2}}^T$ to
\hyperlink{MP}{($\mathrm{MP}$)} exists. Moreover, under the previous condition, strong
existence holds for the nonlinear stochastic differential equation
(\ref{nonlinSDE}) in $[0,T]$, and by \textup{(a)} and \textup{(b)}, one
has $P=\operatorname{law}(\tau,X,\Phi)$. Finally, $\rho^P$ is the unique
solution in $\mathbf{F}_{0,1,T}\cap\mathbf{F}_{0,p,T}$ to the vortex equation
(\ref{eqmild}).
\end{enumerate}
\end{teorema}

The proof of Theorem \ref{teoMP} requires some preliminary facts
about a scalar problem implicitly included in the vectorial
problem \hyperlink{MP}{($\mathrm{MP}$)}.

\subsection{A nonlinear Fokker--Planck equation with external
field associated with the 3d-vortex equation} Recall that the
notation $\tilde{Q}_t$ and $Q^{\circ}_t$ were, respectively,
defined in Definition \ref{margsP} and (\ref{measvect}).

For any $Q\in\mathcal{P}(\mathcal{C}_T)$, we now denote by $\hat{Q}_t$
the sub-probability measure on $\mathbb{R}^3$ defined for scalar
functions by
%
%
\begin{equation}\label{subproba}
\hat{Q}_t(f)=\mathbb{E}^{Q}\bigl(f(X_t) \mathbf{1}_{\{\tau\leq t\}}\bigr),
\end{equation}
where $(\tau,X)$ are the two first marginal of the canonical
process $(\tau,X,\Phi)$ in $\mathcal{C}_T$. Obviously, for $Q\in
\mathcal{P}_b(\mathcal{C}_T)$ we have
\[
\tilde{Q}_t\ll\hat{Q}_t \ll Q^{\circ}_t,
\]
and we shall denote
%
%
\begin{equation}\label{density'}
k^Q_t(x):=\frac{d \tilde{Q}_t}{d \hat{Q}_t}(x).
\end{equation}
Notice that, indeed,
\[
k^Q_t(x)=\frac{E^Q(\Phi_t h(\tau,X_0)\mathbf{1}_{\{\tau\leq t\}}
\vert
X_t=x)}{Q(\tau\leq t\vert X_t=x)} \mathbf{1}_{\{Q(\tau\leq t\vert
X_t=x)>0\}}.
\]
\begin{definicion}
If $Q^{\circ}_t(dx)$ has a density $\rho^Q(t,x)$ with respect to
Lebesgue measure, we shall denote by $\hat{\rho}^Q(t,x)$ the
family of densities of $\hat{Q}_t $.
\end{definicion}

Notice that one has
\[
\tilde{\rho}^Q(t,x)= k^Q_t(x)\hat{\rho}^Q(t,x).
\]
\begin{rem}\label{medibilidad} If $Q\in\mathcal{P}_b(\mathcal
{C}_T)$ is such that $Q_t$ is absolutely continuous for all $t\in
[0,T]$, the existence of a joint measurable version of
$(t,x)\mapsto\rho^Q(t,x)$ is standard by continuity of $X_t$
under $Q^{\circ}_t$. We always work with such a version. Moreover,
there exist measurable versions of $(t,x)\mapsto
\hat{\rho}^Q(t,x)$ and $(t,x)\mapsto\tilde{\rho}^Q(t,x)$. This
can be seen by Lebesgue derivation (see, e.g., Theorem 3.22 in
\cite{Foll}), taking $\delta\to0$ in the quotients
\[
\frac{Q(\tau\leq t , X_t\in B(x,\delta))}{Q(X_t\in B(x,\delta))}
\quad\mbox{and}\quad \frac{E^Q(\Phi_t
h(\tau,X_0)\mathbf{1}_{\{\tau\leq t\}}, X_t\in B(x,\delta
))}{Q(X_t\in
B(x,\delta))}
\]
and using the previous relation between
$\hat{\rho}^Q(t,x)$ and $k^Q$ [here, $B(x,\delta)$ is the open
ball of radius $r$ centered at $x$].
\end{rem}
\begin{lema}\label{eqFP}
Assume that \textup{\hyperlink{MP}{($\mathrm{MP}$)}} has a solution $P\in\mathcal{P}_b(\mathcal
{C}_T)$ such that $P^{\circ}_t$ has a density for each $t\in
[0,T]$. Let $\hat{\rho}:=\hat{\rho}^P$ and $\tilde{\rho}:=
\tilde{\rho}^P$, respectively, denote the densities of $\hat{P}_t$
and $\tilde{P}_t$ and, moreover, assume that (\ref{integKrho})
holds. We have:
\begin{longlist}
\item The couple
$(\hat{\rho},\tilde{\rho})$ satisfies the weak evolution equation
%
%
\begin{eqnarray}\label{weakFP}
&&
\int_{\mathbb{R}^3}f(t,y)\hat{\rho}(t,y)\,dy\nonumber\\
&&\qquad=\int_{\mathbb{R}^3}
f(0,y)\bar
{w}_0(y)\,dy +\int_0^t \int_{\mathbb{R}^3} f(s,y)\bar{\mathbf{g}}(s,y)\,dy
\,ds\nonumber\\[-8pt]\\[-8pt]
&&\qquad\quad{} + \int_0^t \int_{\mathbb{R}^3} \biggl[ \frac{\partial f}{\partial s}(s,y)
+\nu\triangle f(s,y) \nonumber\\
&&\qquad\hspace*{59pt}{} + \mathbf{K}(\tilde{\rho})(s,y)\nabla f(s,y)
\biggr]\hat{\rho}(s,y) \,dy \,ds,\nonumber
\end{eqnarray}
for all $f\in C^{1,2}_b$, where $\bar{w}_0 $ and $\bar{\mathbf{g}}$ were
defined in (\ref{initlaw}).
\item$\hat{\rho}$ is, moreover, a solution of the mild
equation in $[0,T]$,
%
%
\begin{eqnarray}\label{eqmildFP}
\hat{\rho}(t,x) & = & G^{\nu}_t*\bar{w}_0(x)+\int_0^t G^{\nu
}_{t-s}* \bar{\mathbf{g}}(s,\cdot)(x)\,ds\nonumber\\[-8pt]\\[-8pt]
&&{} + \int_0^t\sum_{j=1}^3\int_{\mathbb{R}^3}
\frac{\partial G^{\nu}_{t-s}}
{\partial y_j}
(x-y)\mathbf{K}(k \hat{\rho})_j(s,y)\hat{\rho}(s,y) \,dy \,ds
,\nonumber
\end{eqnarray}
with the multiple integral being absolutely convergent, and where
$k:=k^P$ is the function defined in (\ref{density'}).
\end{longlist}
\end{lema}
\begin{pf} (i)
By the definition of \hyperlink{MP}{($\mathrm{MP}$)} and the fact that
$ \mathbf{1}_{\{\tau\leq t\}}$ is $\mathcal{F}_0$-measurable,
we deduce that the
expectation of the expression
\begin{eqnarray*}
&&
f(t,X_t)\mathbf{1}_{\{t\geq\tau\}}-f (\tau,X_0)\mathbf{1}_{\{t\geq
\tau\}}\\
&&\qquad{}
-\int_0^t \biggl[\frac{\partial f}{\partial s}(s,X_s)+\nu\triangle
f (s,X_s)\,ds + \mathbf{K}(\tilde{\rho})(s,X_s)\nabla
f(s,X_s) \biggr]\mathbf{1}_{\{s\geq\tau\}} \,ds
\end{eqnarray*}
vanishes (see also the beginning of the proof of Lemma
\ref{MPtoweakeq}). Taking expectation and recalling the definition of
$\hat{\rho}$ and $P_0$ [cf. (\ref{subproba}) and (\ref{P0})], we
obtain the desired result applying Fubini's theorem in the time
integral, which is possible since
\[
{\int_{[0,T]\times\mathbb{R}^3}}|\mathbf{K}(\tilde{\rho})(t,x)|\hat
{\rho}(t,x)\,dx
\,dt<\infty,
\]
thanks to condition (\ref{integKrho}).

(ii) Fix $\psi\in\mathcal{D}$ and $t\in[0,T]$ and take in
(\ref{weakFP}) the $C^{1,2}_b$-function $f_t\dvtx[0,t]\times
\mathbb{R}^3\to\mathbb{R}^3$ given by $f_t(s,y)=G^{\nu}_{t-s}*\psi
(y)$ (which
solves the backward heat equation on $[0,t]\times\mathbb{R}^3$ with final
condition $f_t(t,y)=\psi(y)$). By Lemma \ref{estimsw} and
condition~(\ref{integKrho}), it is not hard to check that
\[
\int_0^t\int_{(\mathbb{R}^3)^2} \sum_{j=1}^3 \biggl|\frac{\partial
G^{\nu}_{t-s}}{\partial y_j}(x-y) \biggr||\mathbf{K}(\tilde{\rho})_j(s,y)|
|\psi(x)|\rho(s,y)\,dx \,dy \,ds<\infty.
\]
By Fubini's theorem we easily conclude.
\end{pf}

Consider now a fixed but arbitrary function $k\dvtx[0,T]\times
\mathbb{R}^3\to\mathbb{R}^3$ of class $L^{\infty}([0,T],(L^{\infty
})^3)$, and
formally define an operator $\mathbf{b}^k$ on functions $\eta,\nu\in
\mathcal{M}\mathit{eas}^T$ by
\[
\mathbf{b}^k(\eta,\nu)(t,x)
=\int_0^t\sum_{j=1}^3\int_{\mathbb{R}^3}\frac{\partial
G^{\nu}_{t-s}}{\partial y_j} (x-y)\mathbf{K}(k
\nu)_j(s,y)\eta(s,y) \,dy \,ds.
\]
\begin{rem}\label{FF}
For each $p\in[1,\infty]$ (resp., each $p\in[1,\infty]$ and
$r\geq p$), the mapping $\eta\mapsto k\eta$ is continuous from
$F_{0,p,T}$ to $\mathbf{F}_{0,p,T}$ (resp., from $F_{0,r,(T;p)}$ to
$\mathbf{F}_{0,r,(T;p)}$).
\end{rem}

Write now
\[
\gamma_0(t,x):= G^{\nu}_t*\bar{w}_0(x)+\int_0^t G^{\nu}_{t-s}*
\bar{\mathbf{g}}(s,\cdot)(x)\,ds,
\]
where $\bar{w}_0$ and $\bar{\mathbf{g}}$ were defined in
(\ref{initlaw}). We can state the following properties of the
scalar equation (\ref{eqmildFP}).
\begin{proposicion}\label{contb}
Assume \hyperlink{Hp}{$(\mathrm{H}_1)$} and \hyperlink{Hp}{$(\mathrm{H}_p)$} with $p\in
(\frac{3}{2},3)$, and let $k\in L^{\infty}([0,T],(L^{\infty})^3)$
be a fixed but arbitrary function.
\begin{longlist}
\item For each $r\in[p,\infty)$, we have
\[
\gamma_0\in F_{0,r,(T;p)}\qquad \mbox{with } \tn\gamma_0
\tn_{0,r,(T;p)}\leq C(r,p)\|\bar{w}_0\|_p + T \tn
\bar{\mathbf{g}}\tn_{0,p,T}
\]
for some finite constant $C(r,p)>0$.
\item Suppose that $\frac{3}{2}< p<3, p\leq l <\min\{\frac
{6p}{6-p},3\}$ and $ \frac{3l}{6-l}\leq
l'< \frac{3l}{6-2l}$. Then, there exists a finite constant
$C_0(l,l';p)$ not depending on $T>0$ such that for all
$\eta,\nu\in F_{0,l,(T;p)}$,
\[
\tn\mathbf{b}^k(\eta,\nu) \tn_{0,l',(T;p)} \leq C_0(l,l';p)
T^{1-{3}/({2p})} \tn\eta
\tn_{0,l,(T;p)} \tn\nu\tn_{0,l,(T;p)}.
\]

\item The mild Fokker--Planck
equation with external field (\ref{eqmildFP}) has, at most, one
solution $\hat{\rho}\in F_{0,p,T}$ for each $T>0$.

\item If $\hat{\rho}\in F_{0,p,T}$ is a solution of
(\ref{eqmildFP}), then $\hat{\rho}\in F_{0,r,(T;p)}$ for all $ r
\in[p, \infty)$ with $\tn\hat{\rho} \tn_{0,r,(T;p)}\leq
C(T,p,r,\tn\hat{\rho} \tn_{0,p,T}) < \infty$.

\item We deduce that for all $l\in
[\frac{3p}{3-p},\infty)$, $\mathbf{K}(k \hat{\rho})\in
\mathbf{F}_{1,l,(T;{3p}/({3-p}))}$.
\end{longlist}
\end{proposicion}
\begin{pf} Part (i) follows from Lemma \ref{estimsw} in
a similar way as part (i) of Lem\-ma~\ref{contB}. We notice
that the restriction on $r$ in the latter was needed only to
ensure that the derivative of time integral was convergent, and so
it is not needed here. Thanks to Remark \ref{FF}, part (ii)
is similar to part (ii) of Proposition 3.1 in \cite{F1}.

From the previous parts, equation (\ref{eqmildFP}) admits the
abstract formulation in $F_{0,p,T}$
\[
\hat{\rho}= \gamma_0+ \mathbf{b}^k(\hat{\rho},\hat{\rho}).
\]
Then, the arguments yielding parts (i) of Theorems \ref{exist1} and
\ref{regularity} also provide
the assertions of parts (iii) and (iv), respectively. For
part (v), we notice that from (iv), $k\hat{\rho}\in
\mathbf{F}_{0,r,(T;p)}$ holds for all $r\in[ p,\infty[$. Thus, if we take
$l\geq q:=\frac{3p}{3-p}$ and set
$r:=(\frac{1}{l}+\frac{1}{3})^{-1}$, then one has $r\geq p$, and
so Lemma \ref{contK}(i) implies that
\[
\sup_{t\in
[0,T]}t^{{3/2}({1/p}-{1/r})}\|\mathbf{K}(k\hat
{\rho
})(t,\cdot)\|_l
=\sup_{t\in
[0,T]}t^{{3/2}({1/q}-{1/l})}\|\mathbf{K}(k\hat
{\rho
})(t,\cdot)\|_l
<\infty.
\]
This shows that $\mathbf{K}(k\hat{\rho}) \in\mathbf
{F}_{0,l,(T;q)}$. We conclude
that $\mathbf{K}(k\hat{\rho}) \in\mathbf{F}_{1,l,(T;q)}$, noting
that since
$k\hat{\rho}\in\mathbf{F}_{0,l,(T;p)}$ for all $l\geq q$,
Lemma \ref{contgradK}(i) implies that $\frac{\partial
\mathbf{K}(k\hat{\rho})}{\partial x_k}\in\mathbf{F}_{0,l,(T;p)}$
for all
$k=1,2,3$. In other words,
\begin{eqnarray*}
&&\sup_{t\in
[0,T]}t^{{3/2}({1/p}-{1/l})} \biggl\|\frac
{\partial
\mathbf{K}(k\hat{\rho})(t,\cdot)}{\partial x_k} \biggr\|_l
\\
&&\qquad=\sup
_{t\in
[0,T]}t^{{1/2}+{3/2}({1/q}-{1/l})} \biggl\|
\frac{\partial
\mathbf{K}(k\hat{\rho})(t,\cdot)}{\partial x_k} \biggr\|_l <\infty,
\end{eqnarray*}
which is the required estimate.
\end{pf}

\subsection{Uniqueness in law and pathwise uniqueness}\label{sec42}

We need the following version of Gronwall's lemma:
\begin{lema}\label{Gronesp}
Let $g$ and $k$ be positive functions on $[0,T]$, such that\break
$\int_0^T k(s) \,ds < \infty$, $g$ is bounded, and
\[
g(t)\leq C +\int_0^t g(s) k(s)\,ds \qquad\mbox{for all }t\in[0,T].
\]
Then, we have
\[
g(t)\leq C\exp\int_0^T k(s)\,ds \qquad\mbox{for all }t\in[0,T].
\]
\end{lema}

We are ready to prove parts (a) and (b) in Theorem
\ref{teoMP}.
\begin{pf*}{Proof of Theorem \ref{teoMP}}
Let $P\in\mathcal{P}^T_{b,{3/2}}$ be a solution of \hyperlink{MP}{($\mathrm{MP}$)}.
Since $\rho\in
F_{0,1,T}\cap F_{0,p,T}$, by interpolation we have $\rho\in
F_{0,{3/2},T}$. By Lemma \ref{contK}(i) we deduce
that (\ref{integKrho}) holds.
Moreover, by Lemma \ref{eqFP}(ii), Proposition
\ref{contb}(iv) and Lemma \ref{contgradK}(i), we have
that $\nabla\mathbf{K}(\tilde{\rho})\in F_{0,3,(T;p)}$, and, consequently,
condition (\ref{integgradKrho}) also holds. By Lemma
\ref{MPtoweakeq} we deduce that $\tilde{\rho}$ is a weak solution
of the vortex equation, and, since $k^P_t$ is bounded, we have
$\tilde{\rho}\in\mathbf{F}_{0,p,T}$.

We now need to prove that the latter implies that
$\tilde{\rho}\in\mathbf{F}_{0,p,T}$ is uniquely determined. By Theorem
\ref{exist1}(a) this will follow by checking that
$\tilde{\rho}$ is also mild solution. For fixed $\psi\in(\mathcal
{D})^3$ and $t\in[0,T]$, define $\mathbf{f}_t\dvtx[0,t]\times\mathbb
{R}^3\to
\mathbb{R}^3$
by $\mathbf{f}_t(s,y)=G^{\nu}_{t-s}*\psi(y)$, which is a function of class
$(C^{1,2}_b)^3$ that solves the backward heat equation on
$[0,t]\times\mathbb{R}^3$ with final condition $\mathbf{f}(t,y)=\psi(y)$.
One can
thus take $\mathbf{f}_t$ in the weak vortex equation and, thanks to
conditions (\ref{integKrho}) and (\ref{integgradKrho}),
apply
Fubini's theorem to deduce [since $\psi\in(\mathcal{D})^3$ is
arbitrary] that
\begin{eqnarray*}
&&\tilde{\rho}(t,x)=\mathbf{w}_0(t,x)+
\int_0^t\sum_{j=1}^3\int_{\mathbb{R}^3} \biggl[
\frac{\partial G^{\nu}_{t-s}}
{\partial y_j}
(x-y)[\mathbf{K}(\tilde{\rho})_j(s,y)\tilde{\rho}(s,y)]\\
&&\hspace*{144.1pt}{} +G^{\nu}_{t-s}
(x-y)\biggl[\tilde{\rho}_j(s,y)\,\frac{\partial\mathbf{K}(\tilde{\rho
})}{\partial
y_j}(s,y)\biggr] \biggr]
\,dy \,ds.
\end{eqnarray*}
Since $\tilde{\rho}$ is divergence-free, to see that
$\tilde{\rho}$ solves the mild equation it is enough to justify
an integration by parts of the last term in the previous equation.
We cannot do that at this point since we cannot ensure enough
(Sobolev) regularity of $\tilde{\rho}$. But noting that for
$q=\frac{3p}{3-p}$ one has $1<q^*<\frac{3}{2}$, we see that the
function $\tilde{\rho}=k^P \hat{\rho}$ belongs to $\mathbf
{F}_{0,q^*,T}$ by
interpolation. On the other hand, one has
$G^{\nu}_{t-s} (x-\cdot)\mathbf{K}(\tilde{\rho})(s,\cdot)\in(W^{1,q})^3$
thanks to Proposition \ref{contb}(v). Since by hypothesis,
$\operatorname{div} \tilde{\rho}(s)=0$ in the distribution sense, the fact that
$\tilde{\rho}(s)\in( L^q)^3$ and a density argument allow us to
check
that
\[
\sum_{j=1}^3\int_{\mathbb{R}^3} \tilde{\rho}_j(s,y)\, \frac
{\partial}
{\partial y_j}[G^{\nu}_{t-s}(x-y)\mathbf{K}(\tilde{\rho})(s,y)]\,dy=0
\]
for all
$s\in\ ]0,T]$. Thus, $\mathbf{w}:=\tilde{\rho}$ is the unique
solution of
(\ref{eqmild}) in $\mathbf{F}_{0,p,T}$.

Now, by a standard argument using the semi-martingale
decomposition of the coordinate processes $X^i$ and their products
$X^iX^j$, we obtain that the martingale part of $f(t,X_t)$ in
\hyperlink{MP}{($\mathrm{MP}$)} is given by the stochastic integral $ \sqrt{2\nu}
\int_0^t\nabla f(s,X_s)
\mathbf{1}_{\{s\geq\tau\}} \,dB_s, $ with
respect to
a Brownian motion $B$ defined on some extension of the canonical
space. From this and the previously established uniqueness of
$\tilde{\rho}$, $P$ is the law of a weak solution of the
stochastic differential equation
%
%
\begin{eqnarray}\label{linSDE}
\mbox{\textup{(i)} }&& X_t=X_0+\sqrt{2\nu} \int_0^t \mathbf{1}_{\{s\geq\tau\}
} \,dB_s
+\int_0^t \mathbf{K}(\mathbf{w})
(s,X_s)\mathbf{1}_{\{s\geq\tau\}}\,ds ,\nonumber\\[-8pt]\\[-8pt]
\mbox{\textup{(ii)} }&& \Phi_t=I_3+\int_0^t \nabla\mathbf{K}(\mathbf
{w})(s,X_s)\Phi_s
\mathbf{1}_{\{s\geq\tau\}} \,ds.\nonumber
\end{eqnarray}
Since (\ref{linSDE}) is \textit{linear} in the sense of
McKean, to conclude uniqueness in law it is enough to prove
pathwise uniqueness for it. This is done first for $X$ and then
for $\Phi$, both with help of the estimate on $ \| \nabla
\mathbf{K}(\mathbf{w})(t) \|_{\infty}$ in Theorem \ref
{regularity} and
Gronwall's lemma.
\end{pf*}

\subsection{Pathwise convergence of the mollified processes and strong
existence for small time}

To prove part (c) of Theorem \ref{teoMP}, we shall construct
a strong solution to the nonlinear SDE of part (b) therein.
We shall do so via approximation by solutions to nonlinear SDEs
with regular drift terms $\mathbf{K}^{\varepsilon}(\mathbf
{w}^{\varepsilon})$ and
$\nabla\mathbf{K}^{\varepsilon}(\mathbf{w}^{\varepsilon})$, where
for each
$\varepsilon>0$, $\mathbf{w}^{\varepsilon}\in F_{1,p,T}\cap
F_{0,1,T}$ is
given by Theorem \ref{exist1}. Thus, our
arguments improve the ones developed in \cite{F1} by providing
a pathwise approximation result at an explicit rate. This will be the
key to
carry out the additional improvements on that work in the forthcoming sections.

If $q=\frac{3p}{3-p}$, H\"{o}lder's inequality and the properties
of $\mathbf{K}$ imply that that for all $t\in[0,T]$,
\begin{eqnarray*}
\|\mathbf{K}^{\varepsilon}(\mathbf{w}^{\varepsilon})(t,\cdot)\|
_{\infty}&\leq& C
\|\varphi_{\varepsilon}\|_{q^*}\tn\mathbf{K}(\mathbf
{w}^{\varepsilon})\tn
_{0,q,T}\\
&\leq&
C \|\varphi_{\varepsilon}\|_{q^*}\tn\mathbf{w}^{\varepsilon}\tn_{0,p,T}.
\end{eqnarray*}
Similarly, one has $\|\nabla
\mathbf{K}^{\varepsilon}(\mathbf{w}^{\varepsilon})(t)\|_{\infty
}\leq C \|\nabla
\varphi_{\varepsilon}\|_{q^*}\tn\mathbf{w}^{\varepsilon}\tn
_{0,p,T}$ and
analogous estimates hold for all derivatives. Thus, for each
$\varepsilon>0$, the function $(s,y)\mapsto
\mathbf{K}^{\varepsilon}(\mathbf{w}^{\varepsilon})(s,y)$ is bounded
and continuous
in $y\in\mathbb{R}^3$, and has infinitely many derivatives in $y\in
\mathbb{R}^3$, which are uniformly bounded in $[0,T]\times\mathbb{R}^3$.

We fix now the time interval $[0,T]$ given by Theorem \ref{teoMP}.
It will be useful to consider in what follows the stochastic flow
%
%
\begin{eqnarray}\label{stochflown}
\xi_{s,t}^{\varepsilon}(x)&=&x+\sqrt{2\nu}
(B_t-B_s)\nonumber\\[-8pt]\\[-8pt]
&&{}+\int_s^t
\mathbf{K}^{\varepsilon}(\mathbf{w}^{\varepsilon})(\theta,\xi
^{\varepsilon
}_{s,\theta}(x))\,
d\theta\qquad\mbox{for all } t\in[s,T],\nonumber
\end{eqnarray}
which has a version that is continuously differentiable in $x$ for
all $(s,t)$ thanks to the previously mentioned regularity
properties of $\mathbf{K}^{\varepsilon}(\mathbf{w}^{\varepsilon})$
(cf. Kunita
\cite{Kun}).

We also consider the strong solution of the stochastic
differential equation in $[0,T]$,
%
%
\begin{eqnarray}\label{nonlinSDEreg}
X^{\varepsilon}_t&=&X_0+\sqrt{2\nu} \int_0^t
\mathbf{1}_{\{s\geq\tau\}} \,dB_s+\int_0^t
\mathbf{K}^{\varepsilon}(\mathbf{w}^{\varepsilon})
(s,X^{\varepsilon}_s)\mathbf{1}_{\{s\geq\tau\}}
\,ds,\nonumber\\[-8pt]\\[-8pt]
\Phi^{\varepsilon}_t&=&I_3+\int_0^t \nabla
\mathbf{K}^{\varepsilon}(\mathbf{w}^{\varepsilon})(s,X^{\varepsilon
}_s)\Phi
^{\varepsilon}_s \mathbf{1}_{\{s\geq\tau\}}\,ds,\nonumber
\end{eqnarray}
where $(\tau,X_0)$ is independent of $B$. We denote by
$P^{\varepsilon}$ the joint law of
$(\tau,X^{\varepsilon},\Phi^{\varepsilon})$ and observe that
$P^{\varepsilon}\in\mathcal{P}_b^T$. Since $X^{\varepsilon}_t=X_0$
for all $t\leq\tau$, we have that
\[
X_t^{\varepsilon}=\xi^{\varepsilon}_{\tau,t}(X_0)\mathbf{1}_{\{
t\geq\tau\}
}+X_0\mathbf{1}_{\{t<\tau\}}.\vadjust{\goodbreak}
\]

Denoting by $G^{\varepsilon}(s,x;t,y),(s,x,t,y)\in(\mathbb
{R}_+\times
\mathbb{R}^2)^2, s<t$, the density of $\xi^{\varepsilon}_{s,t}(x)$ (which
is a continuous function of $(s,x,t,y)$, see \cite{Fried}), and
conditioning with respect to $(\tau,X_0)$, we obtain for bounded
and measurable functions $f$ that
\begin{eqnarray*}
E(f(X^{\varepsilon}_t)) &=&
\int_0^t\int_{(\mathbb{R}^3)^2}
f(y) G^{\varepsilon}(s,x;y,t)\,dy P_0(ds,dx)\\
&&{} +\int_t^T\int_{\mathbb{R}
^3}f(x) P_0(ds,dx)\\
&=&\int_{\mathbb{R}^3}f(x)\bar{w}_0(x)\,dx\\
&&{} +\int_0^t\int_{\mathbb{R}^3}
\biggl[\int_{\mathbb{R}^3}f(y)G^{\varepsilon}(s,x;t,y)\,dy \biggr]
\bar
{\mathbf{g}}(s,x) \,dx \,ds \\
&&{}+ \int_t^T\int_{\mathbb{R}^3}f(x)\bar{\mathbf{g}}(s,x)\,dx \,ds.
\end{eqnarray*}
Consequently, $X^{\varepsilon}_t$ has a (bi-measurable) family of
densities that we denote by $\rho^{\varepsilon}$. Observe that one
has $\rho^{\varepsilon}(t)\in L^p$ for all $t\in[0,T]$ from the
assumption on $w_0$ and $\mathbf{g}$ and standard Gaussian bounds for
$G^{\varepsilon}(s,x;t,y)$.

The functions $\hat{\rho}^{\varepsilon}$ and
$\tilde{\rho}^{\varepsilon}$ correspond to the densities of,
respectively, the sub-probability measure and the vectorial
measure
\[
f\mapsto E \bigl[f(\xi^{\varepsilon}_{\tau,t}(X_{\tau}))\mathbf
{1}_{\{t\geq
\tau\}} \bigr]
\]
and
\[
\mathbf{f}\mapsto
E \bigl[\mathbf{f}(\xi^{\varepsilon}_{\tau,t}(X_{\tau})) \nabla_x
\xi^{\varepsilon}_{\tau,t}(X_{\tau}) h(\tau,X_0)\mathbf{1}_{\{
t\geq
\tau\}} \bigr].
\]
They are bi-measurable by similar arguments as in Remark
\ref{medibilidad}, and we have $\hat{\rho}^{\varepsilon}(t)\in
L^p$ and $\tilde{\rho}^{\varepsilon}(t)\in L^p_3$.

The assumptions on $\varphi$ ensure the following estimate
concerning the approximations $\varphi_{\varepsilon}$ of the Dirac
mass (see Lemma 4.4 in Raviart \cite{Rav}):
\begin{lema}\label{aproxco1}
Let $\varphi$ be a cutoff function of order $1$. Then, for all
$v\in W^{1,r}$ and $r\in[1,\infty]$, one has
\[
\|v-\varphi_{\varepsilon}*v\|_r\leq C \varepsilon
\sum_{i=1}^3 \biggl\|\frac{\partial v}{\partial x_i } \biggr\|_r.
\]
\end{lema}

We deduce the following result:
\begin{lema}\label{aproxims}
\textup{(i)} We have $
\tilde{\rho}^{\varepsilon}=\mathbf{w}^{\varepsilon}$ and, consequently,
%
%
\begin{equation}\label{estsrhon}
\sup_{\varepsilon>0} \tn\tilde{\rho}^{\varepsilon}
\tn_{0,p,T}<\infty\quad\mbox{and}\quad \sup_{\varepsilon>0} \tn
\hat{\rho}^{\varepsilon} \tn_{0,p,T} <\infty.
\end{equation}

\textup{(ii)} If $\varphi$ is a cutoff function of order $1$, then we
have that
\[
\sup_{t\in[0,T]} t^{{3}/({2p})-{1}/{2}}
\|\mathbf{w}^{\varepsilon}(t)- \mathbf{w}(t)\|_p\leq C(T)\varepsilon
\]
for some finite constant $C(T)$.
\end{lema}
\begin{pf} (i) Since $E ({\int_0^T }|
\mathbf{K}^{\varepsilon}(\mathbf{w}^{\varepsilon
})(t,X_t^{\varepsilon})|\,dt )
$ and
$E ({\int_0^T} |\nabla
\mathbf{K}^{\varepsilon}(\mathbf{w}^{\varepsilon
})(t,X_t^{\varepsilon})|\,dt )$
are finite, we can follow the lines of Lemma \ref{MPtoweakeq} and use Remark
\ref{P0rem} to see that for all $\mathbf{f}\in(C^{1,2}_b)^3$,
%
%
\begin{eqnarray}\label{weakvareps}
&&
\int_{\mathbb{R}^3}\mathbf{f}(t,y)\tilde{\rho}^{\varepsilon}(t,y)\,dy
\nonumber\\
&&\qquad=
\int_{\mathbb{R}^3} \mathbf{f}(0,y)w_0(y)\,dy + \int_0^t
\int_{\mathbb{R}^3}\mathbf{f}(s,y)\mathbf{g}(s,y)\,dy \,ds \nonumber
\\
&&\qquad\quad{} +\int_0^t
\int_{\mathbb{R}^3} \biggl[\frac{\partial\mathbf{f}}{\partial s}(s,y)
+\nu\triangle\mathbf{f}(s,y)\\
&&\qquad\hspace*{58.6pt}{} + \nabla\mathbf{f}(s,y)
\mathbf{K}^{\varepsilon}(\mathbf{w}^{\varepsilon})(s,y)\nonumber\\
&&\qquad\hspace*{58.6pt}{}+\mathbf
{f}(s,y)\nabla
\mathbf{K}^{\varepsilon}(\mathbf{w}^{\varepsilon})(s,y)
\biggr]\tilde{\rho
}^{\varepsilon}(s,y)
\,dy \,ds.\nonumber
\end{eqnarray}
On the other hand, the regularity properties of the stochastic
flow (\ref{stochflown}) imply that for all $\phi\in\mathcal{D}$
and $\theta\in\ ]0,T]$, the Cauchy problem
%
%
\begin{eqnarray}\label{cauchyprobv}\hspace*{28pt}
&&
\frac{\partial}{\partial s}f(s,y) +\nu\Delta
f (s,y)\nonumber\\
&&\qquad{} +\mathbf{K}^{\varepsilon}(\mathbf{w}^{\varepsilon
})(s,y)\nabla f(s,y)=0,\qquad
(s,y)\in{[0,\theta[}\times\mathbb{R}^3,\\
&&f(\theta,y)=\phi(y)\nonumber
\end{eqnarray}
has a unique solution $f$ that belongs to $C^{1,3}_b([0,\theta]
\times\mathbb{R}^3)$ (see Lemma 4.3 in \cite{F1}). One can thus use the
function $\mathbf{f}=\nabla f$ in (\ref{weakvareps}), and after
simple computations obtain, thanks to the null divergence of $w_0$
and $\mathbf{g}(s,\cdot)$, that
\begin{eqnarray*}
&&\int_{\mathbb{R}^3}\nabla\phi(y)\tilde{\rho}^{(n)}(t,y)\,dy\\
&&\qquad=\int_0^t
\int_{\mathbb{R}^3}\nabla\biggl[ \frac{\partial f}{\partial s}(s,y)
+\nu
\triangle f(s,y) +\mathbf{K}^{\varepsilon}(\mathbf{w}^{\varepsilon
})(s,y)\nabla f
(s,y) \biggr]\\
&&\qquad\quad\hspace*{29pt}{}\times\tilde{\rho}^{(n)}(s,y) \,dy \,ds =0
\end{eqnarray*}
for all $\phi\in\mathcal{D}$. Thus, $\operatorname{div}
\tilde{\rho}^{\varepsilon}(t) =0$, and we can adapt the arguments
of Section \ref{sec42} to conclude that $\tilde{\rho}^{\varepsilon}$
solves the linear mild equation
%
%
\begin{equation}\label{linabs}
\mathbf{v}=\mathbf{w}_0+ \mathbf{B}^{\varepsilon}(\mathbf
{v},\mathbf{w}^{\varepsilon}) ,\qquad \mathbf{v}
\in
\mathbf{F}_{0,p,T}.
\end{equation}
Since uniqueness for (\ref{linabs}) holds (by similar arguments as
for the nonlinear version), and $\mathbf{w}^{\varepsilon}$ also
solves the
equation, we conclude that
$\tilde{\rho}^{\varepsilon}=\mathbf{w}^{\varepsilon}$. The asserted
uniform bound for $\tilde{\rho}^{\varepsilon}$ is thus granted by
Theorem \ref{exist1}. To obtain the uniform bound for
$\hat{\rho}^{\varepsilon}$, we take $L^p$ norm to (\ref{linabs}),
and follow the arguments of the proof of Theorem
\ref{exist1}(i), to get that
\[
\|\tilde{\rho}^{\varepsilon}(t)\|_p \leq\tn\mathbf{w}_0\tn_{0,p,T}+C
\tn\mathbf{w}^{\varepsilon}\tn_{0,p,T}\int_0^t (t-s)^{-{3}/({2p})}
\|\tilde{\rho}^{{\varepsilon}}(s)\|_p \,ds.
\]
The conclusion follows by a similar application of Gronwall's
lemma as therein.

(ii) By an iterative argument as in the proof of Theorem
\ref{exist1}(i), we get that
%
%
\begin{eqnarray}\label{convrhowt}
\|\tilde{\rho}^{{\varepsilon}}(t)-\mathbf{w}(t)\|_p &\leq& C \int_0^t
\alpha(t-s)\|\mathbf{K}^{\varepsilon}(\mathbf{w})(s)-\mathbf
{K}(\mathbf{w})(s)\|_q \,ds\nonumber\\[-8pt]\\[-8pt]
&&{}+C(T) \int_0^t
\|\tilde{\rho}^{\varepsilon}(s)-\mathbf{w}(s)\|_q \,ds,\nonumber
\end{eqnarray}
where $\alpha(s)=\sum_{k=1}^{\tilde{N}(p)}s^{k\theta_0-1}$,
$\theta_0=1-\frac{3}{2p}$ and
$\tilde{N}(p)=\lfloor\theta_0^{-1} \rfloor+1$. Integrating in
time and using Gronwall's lemma, Theorem \ref{regularity}(i)
and Lemma \ref{aproxco1}, we obtain that for all $\theta\in
[0,T]$,
\begin{eqnarray*}
\int_0^{\theta}\|\tilde{\rho}^{\varepsilon}(t)-\mathbf{w}(t)\|_p
\,dt &\leq&
C \int_0^T\int_0^t
\alpha(t-s) \|\mathbf{K}^{\varepsilon}(\mathbf{w})(s)-\mathbf
{K}(\mathbf{w})(s)\|_q \,ds \,dt \\
&\leq& C\varepsilon
\int_0^T \sum_{k=1}^{\tilde{N}(p)}
t^{k(1 -{3}/({2p}))-{1/2}}\,dt = \varepsilon C(T).
\end{eqnarray*}
Substituting the latter in (\ref{convrhowt}), we obtain
\begin{eqnarray*}
\|\tilde{\rho}^{{\varepsilon}}(t)-\mathbf{w}(t)\|_p &\leq&
\varepsilon C(T) + C\int_0^t \alpha(t-s) \|\mathbf{K}^{\varepsilon
}(\mathbf{w}
)(s)-\mathbf{K}(\mathbf{w})(s)\|_q \,ds\\
&\leq& \varepsilon C(T)+ Ct^{{1/2}-{3/(2p)}} \varepsilon,
\end{eqnarray*}
and the conclusion follows.
\end{pf}

The proof of Theorem \ref{teoMP}(c) will be completed by
the following result, which, moreover, establishes the strong
pathwise convergence of the nonlinear processes
$(X^{\varepsilon},\Phi^{\varepsilon})$ as $\varepsilon\to0$. We
are inspired here by ideas introduced in \cite{FM}, but we need a
finer use of analytical properties, as we shall improve the rate
of $\varepsilon^{\delta}$ with $\delta\in(0,1)$, that was
obtained therein for a particular choice of kernel. Further
difficulties also will arise because of the additional (and more
singular) drift term of the ``vortex stretching processes''
$\Phi$, proper to dimension $3$.
\begin{proposicion} Let $\varphi$ be a cutoff of order $1$ and
$K^{\varepsilon}$ be defined in terms of $\varphi$ as before. Then,
as $\varepsilon$ goes to $0$, the family of processes $(X^{\varepsilon
}-X_0,\Phi^{\varepsilon})$, \mbox{$\varepsilon
>0$} is Cauchy in the Banach space of continuous processes
$(Y,\Psi)$ with values in $\mathbb{R}^3\times\mathbb{R}^{3\otimes3}$
with finite norm $E({\sup_{t\in[0,T]}} |Y_t| +|\Psi_t|)$.
Moreover, one has
\[
E \Bigl({\sup_{t\in[0,T]}} |X_t-X_t^{\varepsilon}|+
|\Phi_t-\Phi_t^{\varepsilon}| \Bigr)\leq C(T)\varepsilon,
\]
where $(X,\Phi)$ is a solution of the nonlinear s.d.e.
(\ref{nonlinSDE}).
\end{proposicion}
\begin{pf} We observe that the substraction of $X_0$ is only needed to
avoid a moment-type assumption on $X_0$. Let $\varepsilon>\varepsilon
'>0$. We have
%
%
\begin{eqnarray}\label{estimepseps'}
&&
E \Bigl({\sup_{s\leq t}}|X^{\varepsilon}_s-X^{\varepsilon'}_s| \Bigr)
\nonumber\\
&&\qquad\leq
\int_0^t
E \bigl|\bigl(\mathbf{K}^{\varepsilon}(\mathbf{w}^{\varepsilon
})(s,X^{\varepsilon
}_s)-\mathbf{K}^{\varepsilon'}(\mathbf{w}^{\varepsilon
})(s,X^{\varepsilon}_s)\bigr)
\mathbf{1}_{\{s\geq\tau\}} \bigr|\,ds
\nonumber\\[-8pt]\\[-8pt]
&&\qquad\quad{} + \int_0^t
E \bigl|\bigl(\mathbf{K}^{\varepsilon'}(\mathbf{w}^{\varepsilon
})(s,X^{\varepsilon
}_s)-\mathbf{K}^{\varepsilon'}(\mathbf{w}^{\varepsilon
'})(s,X^{\varepsilon}_s)\bigr)
\mathbf{1}_{\{s\geq\tau\}} \bigr|\,ds
\nonumber\\
&&\qquad\quad{} + \int_0^t
E \bigl|\bigl(\mathbf{K}^{\varepsilon'}(\mathbf{w}^{\varepsilon
'})(s,X^{\varepsilon
}_s)-\mathbf{K}^{\varepsilon'}(\mathbf{w}^{\varepsilon
'})(s,X^{\varepsilon'}_s)\bigr)
\mathbf{1}_{\{s\geq\tau\}} \bigr|\,ds.\nonumber
\end{eqnarray}
The third term on the right-hand side of (\ref{estimepseps'}) is bounded
thanks to Theorem \ref{regularity}(iii) by
\[
C\int_0^t s^{-{1/2}-{3/2}({1/p}-{1/r})}
E \Bigl({\sup_{\theta\leq
s}}|X^{\varepsilon}_{\theta}-X^{\varepsilon'}_{\theta}| \Bigr)\,ds
\]
for any fixed $r\in(3,\frac{3p}{3-p})$. Writing
$q=\frac{3p}{3-p}$ and $q^*$ for its H\"{o}lder conjugate, and using
Lemmas \ref{contK} and \ref{aproxims}(ii), we bound the
second term by
\[
\int_0^T
\|\mathbf{K}^{\varepsilon'}(\mathbf{w}^{\varepsilon})(s)-\mathbf
{K}^{\varepsilon'}(\mathbf{w}
^{\varepsilon'})(s)\|_q
\|\hat{\rho}^{\varepsilon}(s)\|_{q^*}\,ds\leq C(T)\varepsilon.
\]
We have used the fact that ${\sup_{\varepsilon>0}} \tn
\hat{\rho}^{\varepsilon} \tn_{0,q^*,T} <\infty$ by interpolation
since $q^*<\frac{3}{2}<p$. By similar arguments, the first term on
the right-hand side of (\ref{estimepseps'}) can be bounded above by
\[
\int_0^T
\|\mathbf{K}^{\varepsilon'}(\mathbf{w}^{\varepsilon})(s)-\mathbf
{K}^{\varepsilon}(\mathbf{w}
^{\varepsilon})(s)\|_q
\|\hat{\rho}^{\varepsilon}(s)\|_{q^*}\,ds\leq C(T)\varepsilon.
\]
Bringing all together and using Gronwall's lemma we deduce that
%
%
\begin{equation}\label{chauchy1}
E \Bigl({\sup_{s\leq T}}|X^{\varepsilon}_t-X^{\varepsilon'}_t| \Bigr)
\leq C(T) \varepsilon.
\end{equation}

Now, notice that Gronwall's lemma and Theorem \ref{regularity}(iii)
imply that the processes $\Phi_t^{\varepsilon}$ are
bounded in $L^{\infty}( [0,T]\times\Omega,dt\otimes
\mathbb{P})$ uniformly in $\varepsilon$. Therefore, we have
%
%
\begin{eqnarray}\label{estimepseps'grad}
&&
E \Bigl({\sup_{s\leq t}}|\Phi^{\varepsilon}_s-\Phi^{\varepsilon
'}_s| \Bigr) \nonumber\\
&&\qquad\leq
C \int_0^t
E \bigl|\bigl(\nabla\mathbf{K}^{\varepsilon}(\mathbf{w}^{\varepsilon
})(s,X^{\varepsilon
}_s)-\nabla\mathbf{K}^{\varepsilon'}(\mathbf{w}^{\varepsilon})
(s,X^{\varepsilon}_s)\bigr) \mathbf{1}_{\{s\geq\tau\}} \bigr|\,ds
\nonumber\\
&&\qquad\quad{} + C  \int_0^t
E \bigl|\bigl(\nabla\mathbf{K}^{\varepsilon'}(\mathbf{w}^{\varepsilon
})(s,X^{\varepsilon}_s)-\nabla\mathbf{K}^{\varepsilon'}(\mathbf
{w}^{\varepsilon'})
(s,X^{\varepsilon}_s) \bigr)\mathbf{1}_{\{s\geq\tau\}} \bigr|\,ds
\\
&&\qquad\quad{} + C \int_0^t
E \bigl|\bigl(\nabla\mathbf{K}^{\varepsilon'}(\mathbf{w}^{\varepsilon
'})(s,X^{\varepsilon}_s)-\nabla\mathbf{K}^{\varepsilon'}(\mathbf
{w}^{\varepsilon'})
(s,X^{\varepsilon'}_s) \bigr)\mathbf{1}_{\{s\geq\tau\}} \bigr|\,ds
\nonumber\\
&&\qquad\quad{} + C \int_0^t
E \Bigl( {|\nabla\mathbf{K}^{\varepsilon'}(\mathbf
{w}^{\varepsilon'})
(s,X^{\varepsilon'}_s) | \sup_{\theta\leq s}}|\Phi
^{\varepsilon}_{\theta}-\Phi^{\varepsilon'}_{\theta}| \Bigr)\,ds.\nonumber
\end{eqnarray}
By Theorem \ref{regularity}(iii), for fixed $r\in(3,q)$ the
last term in the right-hand side of (\ref{estimepseps'grad}) is bounded by
\[
C\int_0^t s^{-{1/2}-{3/2}({1/p}-{1/r})}
E \Bigl( {\sup_{\theta\leq s}}|\Phi^{\varepsilon}_{\theta}-\Phi
^{\varepsilon'}_{\theta}|
\Bigr)\,ds,
\]
and the third one is by
\[
C\int_0^t s^{-{1/2}-{3/2}({1/p}-{1/r})}
E | X^{\varepsilon}_s-
X^{\varepsilon'}_s | \,ds \leq C(T)\varepsilon,
\]
using also the previous estimates on $E| X^{\varepsilon}_s-
X^{\varepsilon'}_s |$. The first term in
(\ref{estimepseps'grad}) is upper bounded by
%
%
\begin{equation}\label{boundp2}
C\int_0^T \|\hat{\rho}^{\varepsilon}(s)\|_{p^*}
\|\nabla\mathbf{K}^{\varepsilon}(\mathbf{w}^{\varepsilon
})(s)-\nabla\mathbf{K}
^{\varepsilon'}(\mathbf{w}^{\varepsilon})(s)\|_p
\,ds.
\end{equation}
If $p\geq2$, then we have $p^*\leq2$ and so by (\ref{estsrhon})
and interpolation, we deduce that (\ref{boundp2}) is bounded by
\begin{eqnarray*}
&&C\tn\hat{\rho}^{\varepsilon}\tn_{0,p^*,T}\int_0^T \|\nabla
\mathbf{K}
(\varphi_{\varepsilon}*\mathbf{w}^{\varepsilon})(s)-\nabla\mathbf
{K}(\mathbf{w}
^{\varepsilon})\|_p\\
&&\qquad{}+
\|\nabla\mathbf{K}(\mathbf{w}^{\varepsilon})-\nabla\mathbf
{K}(\varphi_{\varepsilon'}*\mathbf{w}
^{\varepsilon})(s)\|_p
\,ds \leq C T \varepsilon.
\end{eqnarray*}
This last inequality is obtained by
Lemmas
\ref{contgradK}(i), \ref{aproxco1}, \ref{aproxims}(i)
and the uniform boundedness
of $(\mathbf{w}^{\varepsilon})_{\varepsilon\geq0}$ in $\mathbf{F}_{1,p,T}$.
If now $\frac{3}{2}<p<2$, then we have $3>p^*>2>p$ and by similar
steps as in the previous case $p\geq2$,
we can upper bound (\ref{boundp2}) by
\begin{eqnarray*}
&&
C \tn\hat{\rho}^{\varepsilon}\tn_{0,p^*,(T;p)} \int_0^T
s^{-{3/2}({1/p}-{1/p^*})}
\|\nabla\mathbf{K}^{\varepsilon}(\mathbf{w}^{\varepsilon
})(s)-\nabla\mathbf{K}
^{\varepsilon'}(\mathbf{w}^{\varepsilon})(s)\|_p
\,ds \\
&&\qquad \leq\varepsilon\sup_{\delta\geq0} \tn
\hat{\rho}^{\delta}\tn_{0,p^*,(T;p)} \int_0^T
s^{-{3/2}({1/p}-{1/p^*})}s^{-{1/2}} \,ds \\
&&\qquad\leq\varepsilon C(T) \sup_{\delta\geq0} \tn
\hat{\rho}^{\delta}\tn_{0,p^*,(T;p)}.
\end{eqnarray*}
We have used here Lemma \ref{aproxco1}, the fact that
$(\mathbf{w}^{\varepsilon})_{\varepsilon\geq0}$ is uniformly
bounded in
$\mathbf{F}_{1,p,T}$ and that
$-\frac{3}{2}(\frac{1}{p}-\frac{1}{p^*})-\frac{1}{2} >-1$ since
$p> \frac{3}{2}$. The fact that the supremum in the previous
estimate is finite, is seen in the same way as part (vi) of
Proposition~\ref{contb}, namely by an iterative argument using the
mild equation (similar as therein) satisfied by
$\hat{\rho}^{\varepsilon}$, starting from the uniform bound in
Lemma \ref{aproxims}(i).

Thus, we have shown that the first term in the right-hand side of
(\ref{estimepseps'grad}) is bounded by a constant times
$\varepsilon$. Let us now tackle the second term in the right-hand side
of~(\ref{estimepseps'grad}). This is bounded by
%
%
\begin{eqnarray}\label{boundp2'}
&&
C\int_0^T \|\hat{\rho}^{\varepsilon}(s)\|_{p^*}
\|\nabla\mathbf{K}^{\varepsilon'}(\mathbf{w}^{\varepsilon
})(s)-\nabla\mathbf{K}
^{\varepsilon'}(\mathbf{w}^{\varepsilon'})(s)\|_p
\,ds \nonumber\\[-8pt]\\[-8pt]
&&\qquad\leq C\int_0^T \|\hat{\rho}^{\varepsilon}(s)\|_{p^*}
\|\mathbf{w}^{\varepsilon}(s)-\mathbf{w}^{\varepsilon'}(s)\|_p
\,ds\nonumber
\end{eqnarray}
thanks to Lemma \ref{contgradK}. By Lemma \ref{aproxims}(ii)
we
can upper bound (\ref{boundp2'}), respectively, by
\[
C \varepsilon\int_0^T s^{{1/2}-{3/(2p)}}\,ds = \varepsilon
C(T)
\]
in the case $p\geq2$, or by
\[
C \varepsilon\int_0^T s^{-{3/2}({1/p}-{1/p^*})}
s^{{1/2}-{3/(2p)}}\,ds=C'(T)\varepsilon
\]
in the case
$p<2$, where the constants are finite since $p>\frac{3}{2}$.

Consequently, we have an estimate of the form
\[
E \Bigl({\sup_{s\leq t}}|\Phi^{\varepsilon}_s-\Phi^{\varepsilon
'}_s| \Bigr)
\leq C \varepsilon+ C\int_0^t s^{-{1/2}-{3/2}(
{1/p}-{1/r})}
E \Bigl( {\sup_{\theta\leq s}}|\Phi^{\varepsilon}_{\theta}-\Phi
^{\varepsilon'}_{\theta}|
\Bigr)\,ds
\]
for each fixed $r \in(3,q)$, and Gronwall's lemma yields
%
%
\begin{equation}\label{chauchy2}
E \Bigl({\sup_{s\leq t}}|\Phi^{\varepsilon}_s-\Phi^{\varepsilon
'}_s| \Bigr) \leq
C(T) \varepsilon
\end{equation}
for all $\varepsilon\geq\varepsilon'>0$.

Estimates (\ref{chauchy1}) and (\ref{chauchy2}) thus show that
$(X^{\varepsilon}-X_0,\Phi^{\varepsilon})$ is a
Cauchy sequence in the Banach space of continuous processes
$(Y,\Psi)$ with values in $\mathbb{R}^3\times\mathbb{R}^{3\otimes3}$
and finite norm $E({\sup_{t\in[0,T]}} |Y_t| +|\Psi_t|)$.
Write the limit in the form $(X-X_0,\Phi)$, for a continuous
process $(X,\Phi)$ and define $\mathcal{E}_t^1$ and $\mathcal
{E}_t^2$ by
the relations
%
%
\begin{eqnarray}\label{linsdeE}
X_t&=&X_0+\sqrt{2\nu} \int_0^t \mathbf{1}_{\{s\geq
\tau\}}\,
dB_s+\int_0^t \mathbf{K}(\mathbf{w}) (s,X_s)\mathbf{1}_{\{s\geq
\tau\}} \,ds+\mathcal
{E}_t^1,\nonumber\\[-8pt]\\[-8pt]
\Phi_t&=&I_3+\int_0^t \nabla
\mathbf{K}(\mathbf{w})(s,X_s)\Phi_s\mathbf{1}_{\{s\geq\tau\}}
\,ds+\mathcal{E}_t^2.\nonumber
\end{eqnarray}
Comparing $(X,\Phi)$ and $(X^{\varepsilon},\Phi^{\varepsilon})$,
and using similar estimates as so far in this proof, but with $0$
instead of $\varepsilon'$ (and $\mathbf{w}$ instead of
$\mathbf{w}^{\varepsilon'}$), we get that $(X,\Phi)$ satisfies
(\ref{linsdeE}) with $\mathcal{E}_t^i=0$, $i=1,2$. Since that is a
linear s.d.e. (in McKean's sense), the proof that $(X,\Phi)$ is
the asserted nonlinear process will be achieved by checking that
for all bounded Lipschitz function $\mathbf{f}\dvtx\mathbb{R}^3\to
\mathbb{R}^3$, one has
\[
E\bigl(\mathbf{f}(X_t)\Phi_t h(\tau,X_0)\mathbf{1}_{\{s\geq\tau\}
}\bigr)=\int_{\mathbb{R}
^3}\mathbf{f}(x)\mathbf{w}(t,x)\,dx.
\]
The latter follows from the facts that
\[
E\bigl(\mathbf{f}(X_t^{\varepsilon})\Phi_t^{\varepsilon} h(\tau
,X_0)\mathbf{1}
_{\{s\geq\tau\}}\bigr)=\int_{\mathbb{R}^3}\mathbf{f}(x)\mathbf
{w}^{\varepsilon}(t,x)\,dx
\]
and
%
%
\begin{eqnarray}\label{estimflip}
&&
\bigl|E\bigl(\mathbf{f}(X_t)\Phi_t h(\tau
,X_0)\mathbf{1}_{\{s\geq
\tau\}}\bigr)-E\bigl(\mathbf{f}(X_t^{\varepsilon})\Phi_t^{\varepsilon}
h(\tau,X_0)\mathbf{1}_{\{s\geq\tau\}}\bigr)\bigr| \nonumber\\
&&\qquad\leq
\bigl(\|\Phi\|_{L^{\infty}([0,T]\times
\Omega)}+1\bigr)\|h\|_{\infty}\|\mathbf{f}\|_{\mathrm{Lip}}
E(|X_t-X_t^{\varepsilon}|+|\Phi_t-\Phi_t^{\varepsilon}|)\\
&&\qquad\leq C
\|\mathbf{f}\|_{\mathrm{Lip}} \varepsilon.\nonumber
\end{eqnarray}
\upqed\end{pf}
\begin{rem}
(a)
By Lemma \ref{aproxims}(i), the process
$(X^{\varepsilon},\Phi^{\varepsilon})$ defined in
(\ref{nonlinSDEreg}) is a solution in $[0,T]$ of the nonlinear
s.d.e.:
%
%
\begin{eqnarray}\label{nonlinSDEreg'}
\mbox{(i) }&& X_t^{\varepsilon}=X_0+\sqrt{2\nu} \int_0^t \mathbf
{1}_{\{s\geq
\tau\}} \,dB_s +\int_0^t
\mathbf{K}^{\varepsilon}(\tilde{\rho}^{\varepsilon})
(s,X_s^{\varepsilon})\mathbf{1}_{\{s\geq\tau\}}\,ds ,\nonumber\\
\mbox{(ii) }&& \Phi_t^{\varepsilon}=I_3+\int_0^t \nabla
\mathbf{K}^{\varepsilon}(\tilde{\rho}^{\varepsilon
})(s,X_s^{\varepsilon
})\Phi_s^{\varepsilon}
\mathbf{1}_{\{s\geq\tau\}} \,ds \quad\mbox{and}\nonumber\\[-8pt]\\[-8pt]
\mbox{(iii) }&& \mbox{the law $P^{\varepsilon}$ of
$(\tau,X^{\varepsilon},\Phi^{\varepsilon})$ belongs to }\mathcal
{P}_{b,{3/2}}^T \quad\mbox{and}\nonumber\\
&&\tilde{P}_t^{\varepsilon}(dx)=\tilde{\rho}^{\varepsilon}(t,x)\,dx.\nonumber
\end{eqnarray}

{\smallskipamount=0pt
\begin{longlist}[(b)]
\item[(b)]
It is also possible to associate a unique pathwise solution of
(\ref{nonlinSDE}) with any solution $\mathbf{w}\in\mathbf
{F}_{0,p,T}\cap
\mathbf{F}_{0,1,T}$ of the mild vortex equation (i.e., not necessarily the
one given by Theorem \ref{exist1}). This can be done by an
approximation argument similar to the previous one, but
considering linear processes in the sense of McKean [with drift
terms $\mathbf{K}^{\varepsilon}(\mathbf{w})$ and $\nabla\mathbf
{K}^{\varepsilon}(\mathbf{w})$]
instead of the processes (\ref{nonlinSDEreg}).
\item[(c)]
Denote now by $\mathcal{W}_T$ the Wasserstein
distance in $\mathcal{P}(\mathcal{C}_T)$ associated with the metric in
$\mathcal{C}_T:=[0,T]\times C([0,T],\mathbb{R}^3\times\mathbb
{R}^{3\otimes3})$
\begin{eqnarray*}
&&d((\theta,y,\psi),(\eta,x,\phi)) \\
&&\qquad = |\theta-\eta|+\sup_{t\in[0,T]}
\bigl(\min\{|x(t)-y(t)|,1\}
+ \min\{|\psi(t)-\phi(t)|,1\}\bigr).
\end{eqnarray*}
Then, the previous proof states that
\[
\mathcal{W}_T(P^{\varepsilon},P)\leq C(T) \varepsilon,
\]
where $P$ is the law of the nonlinear process (\ref{nonlinSDE}).
\item[(d)] By the regularity results of Section \ref{sec3}, one can
prove in
a similar way as in Corollary 4.3 of \cite{F1} that the
stochastic flow (\ref{3linearflow}) is of class $C^1$, in spite of
the fact that $\mathbf{u}$ and $\nabla\mathbf{u}$ are singular at
$t=0$. Thus,
identity (\ref{formstocflow}) holds.
\end{longlist}}
\end{rem}

\section{The stochastic vortex method}\label{sec5}

We first consider a McKean--Vlasov model with mollified
interaction and cutoff. This extends the model studied in
\cite{F1} to the present situation involving random space--time
births.

Denote by $M_{\varepsilon}$ the sup-norm of $K_{\varepsilon}$
on $\mathbb{R}^3$
and by $L_{\varepsilon}$ a Lipschitz constant for it, which,
respectively, behave like ${1\over
\varepsilon^3}$ and ${1\over\varepsilon^4}$ when $\varepsilon\ll 1$.
Notice that $\operatorname{div} K_{\varepsilon}=(\operatorname{div}
K) *\varphi_{\varepsilon}=0$.

For $R>0$, we denote by $\chi_R\dvtx\mathbb{R}^{3\otimes3}\to\mathbb
{R}^{3\otimes
3}$ a Lipschitz continuous truncation function such that
$|\chi_R(\phi)|\leq R$. We may and shall assume that $\chi_R$ has
Lipschitz constant less than or equal to $1$.

Consider now a filtered probability space endowed with an adapted
standard three-dimensional Brownian motion $B$ and with a
$[0,T]\times\mathbb{R}^3$-valued random variable $(\tau,X_0)$
independent of $B$ and with law $P_0$.
\begin{teorema}\label{3teoMKV}
There is existence and uniqueness (pathwise and in law) for the
nonlinear process with random space--time births, nonlinear in the
sense of McKean
%
%
\begin{eqnarray}\label{procnonlin}
X^ {\varepsilon, R}_t&=&X_0+\sqrt{2\nu} \int_0^t
\mathbf{1}_{\{s\geq\tau\}} \,dB_s +\int_0^t \mathbf{u}^ {\varepsilon
, R}
(s,X^ {\varepsilon,
R}_s)\mathbf{1}_{\{s\geq\tau\}}\,ds\nonumber\\[-8pt]\\[-8pt]
\Phi^ {\varepsilon, R}_t&=&I_3+\int_0^t \nabla
\mathbf{u}^ {\varepsilon, R}(s,X^ {\varepsilon, R}_s)\chi_R(\Phi^
{\varepsilon, R}_s) \mathbf{1}_{\{s\geq\tau\}} \,ds\nonumber
\end{eqnarray}
with
%
%
\begin{equation}\label{nonlinearite}
\mathbf{u}^ {\varepsilon, R}(s,x)=E \bigl[K_{\varepsilon}(x-X^
{\varepsilon, R}_s)\wedge\chi_R(\Phi^ {\varepsilon, R}_s)
h(\tau,X_0)\mathbf{1}_{\{s\geq\tau\}} \bigr].
\end{equation}
\end{teorema}

The proof is based in the classic contraction argument of
Sznitmann \cite{Szn} and is not hard to obtain by combining
elements of Theorems 5.1 in \cite{F1} and Theorem~3.1 in
\cite{FM}.

Consider next a probability space endowed with a sequence
$(B^i)_{i\in\mathbb{N}}$ of independent three-dimensional Brownian motions,
and a sequence of independent random variables
$(\tau^i,X_0^i)_{i\in\mathbb{N}}$ with law $P_0$ and independent of the
Brownian motions. For each $n\in\mathbb{N}$ and $R,\varepsilon>0$, we
define the following system of interacting particles:
%
%
\begin{eqnarray}\label{lesystemepsR}
X^{i, \varepsilon, R,n}_t & = & X^i_0+\sqrt{2\nu} \int_0^t\mathbf
{1}_{\{s\geq
\tau^i\}} \,dB^i_s \nonumber\\
&&{} + \int_0^t \frac{1}{n}\sum_{j\not=i}
K_{\varepsilon} (X^ {i, \varepsilon,
R,n}_s-X^{j,\varepsilon,R,n}_s)\nonumber\\
&&\hspace*{51.5pt}{}\wedge
\chi_R(\Phi_s^{j,\varepsilon,R,n}) h(\tau^j,X_0^j)\mathbf{1}_{\{
s\geq\tau
^i,\tau^j\}} \,ds,\nonumber\\[-8pt]\\[-8pt]
\Phi^{i ,\varepsilon, R,n}_t & = & I_3 +\int_0^t
\frac{1}{n}\sum_{j\not=i} [\nabla K_{\varepsilon} (X^ {i,
\varepsilon, R,n}_s-X^{j,\varepsilon,R,n}_s)\nonumber\\
&&\hspace*{64.3pt}{} \wedge
\chi_R(\Phi_s^{j,\varepsilon,R,n}) h(\tau^j,X_0^j) ]\nonumber\\
&&\hspace*{58.6pt}{}\times\chi_R(\Phi^
{i, \varepsilon, R,n}_s)
\mathbf{1}_{\{s\geq\tau^i,\tau^j\}}\,ds,\nonumber
\end{eqnarray}
for $i=1,\ldots, n$, and with $\nabla K(y)\wedge z=\nabla_y
(K(y)\wedge z)$ for $y,z\in\mathbb{R}^3,y\not=0$. Pathwise existence
and uniqueness can be proved by adapting standard arguments,
thanks to the Lipschitz continuity of the coefficients.

In
the same probability space, we also consider the sequence
%
%
\begin{eqnarray}\label{procnonlincopies}
X^{i,\varepsilon,R}_t&=&X_0^i+\sqrt{2\nu} \int_0
\mathbf{1}_{\{s\geq\tau^i\}}\,d B^i_s+\int_0^t \mathbf{u}^
{\varepsilon, R}
(s,X^{i,\varepsilon,R}_s)\mathbf{1}_{\{s\geq\tau^i\}}\,ds,\nonumber\\[-8pt]\\[-8pt]
\Phi^{i,\varepsilon,R}_t&=&I_3+\int_0^t \nabla
\mathbf{u}^ {\varepsilon, R}(s,X^{i,\varepsilon,R}_s)\chi_R(\Phi
^{i,\varepsilon,R}_s)\mathbf{1}_{\{s\geq\tau^i\}} \,ds,\qquad
i\in\mathbb{N},\nonumber
\end{eqnarray}
of independent copies of (\ref{procnonlin}). Their common law in
$\mathcal{C}_T$ is denoted by $P^{\varepsilon, R}$, and we write
$\bar{h}:= \|w_0\|_1+\|\mathbf{g}\|_{1,T}$. Recall that $\chi_R$ is a
Lipschitz-continuous function, bounded by $R>0$ and with Lipschitz
constant less than or equal to $1$. It is not hard to adapt the
proof of Theorem 5.2 in \cite{F1} to get the following:
\begin{teorema}\label{propchaosepsR}
For $\varepsilon>0$ sufficiently small and all $R>0$, we have
%
%
\begin{equation}\label{convi}\qquad
\mathbb{E} \Bigl[\sup_{t\in[0,T]} \{
|X^{i,\varepsilon,R,n}_t-X^{i,\varepsilon,R}_t|+
|\Phi^{i,\varepsilon,R,n}_t-\Phi^{i,\varepsilon,R}_t|
\} \Bigr]\leq\frac{1}{\sqrt{n}} C(\varepsilon,R,\bar{h},T)
\end{equation}
for all $i\leq n$, where
\[
C(\varepsilon,R,\bar{h},T)= C_1\varepsilon(1+R \bar{h}
T)(R\bar{h}T)\exp\{C_2\varepsilon^{-9}\bar{h}T(R+1)(\bar{h}+RT)\}
\]
for some positive constants $C_1,C_2$ independent of $R$,
$\varepsilon$, $T$ and $\bar{h}$.
\end{teorema}

Let us now make the assumptions of Theorem \ref{exist1}, and
consider, in the corresponding time interval $[0,T]$, independent
copies $(X^{i,\varepsilon},\Phi^{i,\varepsilon})$ and
$(X,\Phi^i)$ of the processes (\ref{nonlinSDE}) and
(\ref{nonlinSDEreg'}) constructed on the given data
$(X_0^i,\tau^i,B^i)$, $i\in\mathbb{N}$.

Recall again that the uniform bound of Theorem \ref{regularity}(iii)
and Gronwall's lemma imply that the processes
$\Phi^{\varepsilon}$ are uniformly bounded, say
\[
{\sup_{t\in[0,T],\varepsilon\geq0, \omega\in
\Omega}}|\Phi_t^{\varepsilon}(\omega)| \leq R_{\circ}(T,\mathbf{w}_0)
\]
for some finite positive constant $R_{\circ}(T,\mathbf{w}_0)$. Thus, for
any $R\geq R_{\circ}$, one has for all $t\in[0,T]$ that
\[
(X^{i,\varepsilon}_t,\Phi^{i,\varepsilon}_t)=(X^{i,\varepsilon
}_t,\chi_R(\Phi^{i,\varepsilon}_t)).
\]
Consequently, $(X^{i,\varepsilon},\Phi^{i,\varepsilon})$ is a
pathwise solution in $[0,T]$ of (\ref{procnonlincopies}),
and so we conclude that
\[
(X^{i,\varepsilon},\Phi^{i,\varepsilon})=(X^{i,\varepsilon,R},\Phi
^{i,\varepsilon,R})
\]
almost surely. Bringing it all together, we obtain the following
pathwise approximation result:
\begin{teorema}\label{convortmet} Assume that \hyperlink{Hp}{$(\mathrm{H}_1)$} and
\hyperlink{Hp}{$(\mathrm{H}_p)$} hold with $p\in
(\frac{3}{2},3)$ and that the hypothesis of Theorem
\ref{exist1}\textup{(i)}
is satisfied. Let $K_{\varepsilon}$ be defined as in
(\ref{biotsavop}), with $\varphi$ a cutoff function of order $1$
and write $\bar{h}=\|w_0\|_1+\|\mathbf{g}\|_{1,T}$. Let, furthermore,
$R\geq
R_{\circ}(T,\mathbf{w}_0)$ and
\[
\varepsilon_n=(c_{\alpha} \ln n)^{-{1/9}}
\]
with
\[
0<c_{\alpha}<\alpha\bigl(C_2 \bar{h}T(R+1)(\bar{h}+RT)\bigr)^{-1}
\]
for some alpha $\alpha\in(0,\frac{1}{2})$. Then, we have for all
$i\leq n$,
%
%
\begin{eqnarray}\label{convi'}
&&\mathbb{E} \Bigl[\sup_{t\in[0,T]} \{
|X^{i,\varepsilon_n,R,n}_t-X^i_t|+
|\Phi^{i,\varepsilon_n,R,n}_t-\Phi^i_t| \} \Bigr] \nonumber\\[-8pt]\\[-8pt]
&&\qquad\leq
C(T,w_0,\mathbf{g},\alpha) \biggl[\frac{1}{n^{{1/2}-\alpha
}(\ln
n)^{{1/9}}}+\frac{1}{(\ln n)^{{1/9}}} \biggr],\nonumber
\end{eqnarray}
where $(X,\Phi)$ is the unique pathwise solution of
(\ref{nonlinSDE}), and the constant $C(T,w_0,\break\mathbf{g},\alpha)$ depends
on the data $w_0$ and $\mathbf{g}$ only through the quantities $\|w_0\|_p,
\tn\mathbf{g}\tn_{0,p,T}$ and $\|w_0\|_1+\|\mathbf{g}\|_{1,T}$.
\end{teorema}
\begin{rem}
(i)
The rate at which the second term in the right-hand side of
(\ref{convi'}) goes to $0$ is exactly that of
$\varepsilon=\varepsilon_n$. The logarithmic order of latter was
needed to make the upper bound in Theorem \ref{propchaosepsR} go
to $0$ with $n$, which then happens at an algebraic rate. The
global rate is, therefore, conditioned by the
techniques used in the proof of Theorem
\ref{propchaosepsR} (see \cite{F1} for details).
Under additional regularity assumptions, it
is possible by analytic arguments to slightly improve the
convergence rate (see the discussion at the end). An attempt for a
more substantial improvement should, however, exploit specific
features of the interaction at the level of the particle systems.

(ii)
The previous result implies as usual that $\mathcal{W}_T
(\operatorname{law}(X^{i,\varepsilon,R,n},\Phi^{i,\varepsilon,R,n}),P )$
goes to $0$ at least that fast, and that (with the obviously
extended meaning of $\mathcal{W}_T$)
\[
\mathcal{W}_T (\operatorname{law} ((X^{1,\varepsilon,R,n},\Phi
^{1,\varepsilon,R,n}),\ldots,(X^{k,\varepsilon,R,n},\Phi
^{k,\varepsilon,R,n}) )
,P^{\otimes k} )\leq k \delta_n,
\]
where $\delta_n$ stands
for the quantity in the right-hand side of (\ref{convi'}).
\end{rem}

We deduce the convergence at the level of empirical processes:
\begin{corolario}\label{convempproc} Under the assumptions of Theorem
\ref{convortmet}, the family\break
$(\tilde{\mu}^{n,\varepsilon_n,R}_t)_{0\leq t \leq T}$ of
$\mathbb{R}^3$-weighted empirical measures on $\mathbb{R}^3$
\[
\tilde{\mu}^{n,\varepsilon_n,R}_t:= \frac{1}{n}\sum_{i=1}^n
\delta_{X^{i,\varepsilon_n,R,n}_t}
\cdot(\chi_R(\Phi_t^{i,\varepsilon_n,R,n})h_0(\tau,X_0^i)
)\mathbf{1}_{\{t\geq r\}}
\]
converges in probability to $(\mathbf{w}(t,x)\,dx)_{0\leq t
\leq T}$ in the space
$C([0,T],\mathcal{M}_3(\mathbb{R}^3))$, where $\mathcal{M}_3(\mathbb
{R}^3)$ denotes
the space of finite $\mathbb{R}^3$-valued measures on $\mathbb{R}^3$ endowed
with the weak topology. Moreover, we have
\begin{eqnarray*}
&&\sup_{t\in[0,T],\|\mathbf{f}\|_{\mathrm{Lip}}\leq1}E |\langle\tilde
{\mu}^{n,\varepsilon_n,R}_t- \mathbf{w}(t),\mathbf{f}\rangle
| \\
&&\qquad\leq C
\biggl[\frac{1}{\sqrt{n}}+\frac{1}{n^{{1/2}-\alpha}(\ln
n)^{{1/9}}}+\frac{1}{(\ln n)^{{1/9}}} \biggr],
\end{eqnarray*}
where $\| \mathbf{f}\|_{\mathrm{Lip}}$ is the usual norm in the space of
bounded Lipshitz continuous functions $\mathbf{f}\dvtx\mathbb{R}^3 \to
\mathbb{R}^3$.
\end{corolario}
\begin{pf} It is enough to prove the bound for Lipshitz bounded
functions. For such a function
$\mathbf{f}\dvtx\mathbb{R}^3\to\mathbb{R}^3$, it holds
that
%
%
\begin{eqnarray}\label{estimatesf}
&&
|\langle\tilde{\mu}^{n,\varepsilon_n,R}_t,\mathbf{f}\rangle
-\langle\mathbf{w}(t),\mathbf{f}\rangle| \nonumber\\
&&\qquad
\leq\Biggl|\langle\tilde{\mu}^{n,\varepsilon_n,R}_t,\mathbf
{f}\rangle
-\frac{1}{n}\sum_{i=1}^n \mathbf{f}( X_t^{i,\varepsilon
_n,R})\wedge
(\chi_R(\Phi_t^{i,\varepsilon_n,R}))h(\tau,X_0^i)\mathbf{1}_{\{
\tau\geq t\}
} \Biggr|\nonumber\\
&&\qquad\quad{} + \Biggl|\frac{1}{n}\sum_{i=1}^n \mathbf{f}( X_t^{i,\varepsilon
_n,R})\wedge
(\chi_R(\Phi_t^{i,\varepsilon_n,R}))h(\tau,X_0^i)\mathbf{1}_{\{
\tau\geq t\}
}\\
&&\qquad\quad\hspace*{14.5pt}{}-\int_{\mathcal{C}_T} \mathbf{f}( y(t))\wedge
\chi_R(\phi(t))h(\theta,x(0))P^{\varepsilon_n,R}(d\theta,dy,d\phi
) \Biggr|\nonumber\\
&&\qquad\quad{} +|\langle\mathbf{w}^{\varepsilon_n}(t)- \mathbf{w}(t),\mathbf
{f}\rangle|\nonumber
\end{eqnarray}
with
$P^{\varepsilon_n,R}=P^{\varepsilon_n}=\operatorname{law}(\tau,X^{i,\varepsilon
_n,R},\Phi^{i,\varepsilon_n,R})$.
The independence of the processes
$(\tau^i,X^{i,\varepsilon_n,R},\Phi^{i,\varepsilon_n,R})$, $i\in
\mathbb{N}$, and the definition of $h$ imply that the expectation of the
second term in the right-hand side of\vspace*{-2pt} (\ref{estimatesf}) is
bounded by $\frac{1}{\sqrt{n}}2\|\mathbf{f}\|_{\mathrm{Lip}}R \bar{h}$, where
$\bar{h}=(\|w_0\|_1+\|\mathbf{g}\|_{1,T})$. We use the latter and estimate
in Theorem \ref{propchaosepsR} to bound the first term, and get
that
\begin{eqnarray*}
&&
E |\langle\tilde{\mu}^{n,\varepsilon_n,R}_t- \mathbf
{w}(t),\mathbf
{f}\rangle
| \\
&&\qquad\leq \|\mathbf{f}\|_{\mathrm{Lip}}( R+1) \bar{h}\frac{1}{\sqrt{n}}
C(\varepsilon_n,R,\bar{h},T)
\\
&&\qquad\quad{} +\frac{2 \|\mathbf{f}\|_{\mathrm{Lip}} R \bar{h} }{\sqrt{n}} + |\langle
\mathbf{w}^{\varepsilon_n}-\mathbf{w}(t),\mathbf{f}\rangle|.
\end{eqnarray*}
The last term being equal to the first term in
(\ref{estimflip}), the conclusion follows.
\end{pf}
\begin{rem} In the
case $\mathbf{g}=0$, Philipowski \cite{Phi} obtained a similar
approximation result of the vorticity field, for a simpler
particle system, under the additional assumption that the test
function $\mathbf{f}$ belongs to $L^{p^*}$.
\end{rem}

Finally, we establish an approximation result with convergence
rate for the velocity field. To that end, we need to strengthen
the already shown convergence of $\mathbf{w}^{\varepsilon}$ to
$\mathbf{w}$. We
will need the following:
\begin{lema}\label{convgradw} For each $\tilde{p}\in(\frac
{3}{2},p)$, there is a
constant $C(T,\tilde{p})$ such that
\[
\sup_{t\in[0,T]}t^{{3}/({2\tilde{p}})}\|\nabla
\mathbf{w}^{\varepsilon}(t)-\nabla\mathbf{w}(t)\|_{\tilde{p}}\leq
C(T,\tilde{p})
\varepsilon.
\]
\end{lema}
\begin{pf} We need $\tilde{p}\in(\frac{3}{2},3)$ in order to
dispose from a integrable (in time) bound for $\| D^2
\mathbf{K}^{\varepsilon}(\mathbf{w}^{\varepsilon})(t)\|_{
{3\tilde
{p}}/({3-\tilde{p}})}$,
which we do not have for $\tilde{p}=p$. Indeed, for any
$\tilde{p}$ in that interval we have
$\tilde{q}:=\frac{3\tilde{p}}{3-\tilde{p}}\in(3,\frac{3p}{3-p})$,
and so by Theorem \ref{regularity}(i) and Lemma
\ref{contgradK} we have for $k,j,i=1,2,3$ that
%
%
\begin{eqnarray}\label{intboundqtilde}
&&\sup_{t \in[0,T],\varepsilon\geq
0}t^{{3/2}({1/p}-{1}/{\tilde{q}})} \biggl\|
\frac{\partial
\mathbf{K}^{\varepsilon}(\mathbf{w}^{\varepsilon})_j}{\partial
x_i} \biggr\|_{\tilde{q}}\nonumber\\[-8pt]\\[-8pt]
&&\qquad{}+ \sup_{t \in[0,T],\varepsilon\geq
0}t^{{1/2}+{3/2}({1/p}-{1/\tilde{q}})}
\biggl\| \frac{\partial^2
\mathbf{K}^{\varepsilon}(\mathbf{w}^{\varepsilon})_j}{\partial
x_i\,\partial
x_k} \biggr\|_{\tilde{q}}<\infty\nonumber
\end{eqnarray}
with
$-\frac{1}{2}-\frac{3}{2}(\frac{1}{p}-\frac{1}{\tilde
{q}})=-1+\frac{3}{2}(\frac{1}{\tilde{p}}-\frac{1}{p})
>-1$. Let us now check that one has
%
%
\begin{equation}\label{unifbound11}
\sup_{\varepsilon\geq0}\tn
\mathbf{w}^{\varepsilon}\|_{1,\tilde{p},T}<\infty.
\end{equation}
This is not immediate, since $T>0$ given by Theorem \ref{exist1}
was determined by the norm of $\mathbf{w}_0$ and of the operator
$\mathbf{B}^{\varepsilon}$ in the spaces corresponding to the parameter
$p>\tilde{p}$. We will prove (\ref{unifbound11}) using continuity
properties of the operators $\mathbf{B}^{\varepsilon}$. It follows from
Proposition 3.1(iii) in \cite{F1} that for $\frac{3}{2}\leq
r<3$ and $ \frac{3r}{6-r}\leq r' \leq r$, one has
%
%
\begin{equation}\label{contBext}
\sup_{\varepsilon\geq0}\tn\mathbf{B}^{\varepsilon}(\mathbf
{v},\mathbf{v}) \tn_{1,r',T}
\leq C_{r,r'} (T) (\tn\mathbf{v}\tn_{1,r,T})^2
\end{equation}
for some finite constant
$C_{r,r'} (T)$. From this, we deduce that $\mathbf{w}^{\varepsilon}
\in
\mathbf{F}_{1,\tilde{p},T}$, with a uniform (in $\varepsilon$)
bound, by
the following iterative procedure. Define a real sequence by
$r_0=\tilde{p}$, $r_{n+1}=\frac{6r_n}{3+r_n}$, and notice that it
is increasingly convergent to $3$. We can thus take $N\in\mathbb{N}$
such that $r_N< p\leq r_{N+1}$. The function $s\mapsto
\frac{3s}{6-s}$ being increasing on $[0,6]$, we then have
$\frac{3p}{6-p}\leq\frac{3r_{N+1}}{6-r_{N+1}}=r_N$. By
(\ref{contBext}) with $r=p$ and $r'=r_N$, we see that
$\mathbf{B}^{\varepsilon}(\mathbf{w}^{\varepsilon},\mathbf
{w}^{\varepsilon})\in\mathbf{F}_{1,
r_N,T}$, and since also $\mathbf{w}_0 \in\mathbf{F}_{1,r_N,T}$ holds
by Lemma
\ref{contB}(i) (taking $r_N$ in the place of $p$ and $r$
therein), we get that $\mathbf{w}^{\varepsilon}\in\mathbf{F}_{1,
r_N,T}$, with a
bound in that space that is uniform in $\varepsilon$.
We repeat the previous arguments with $r=r_N$ and
$r'=\frac{3r_N}{6-r_N}=r_{N-1}$ and get that $\mathbf{w}^{\varepsilon
}\in
\mathbf{F}_{1, r_{N-1},T}$, with a bound that is a uniform in
$\varepsilon$. Continuing $N-1$ times this scheme we get
(\ref{unifbound11}).

We now take derivatives in the mild vortex equation with
$\varepsilon\geq0$ (as justified in the proof of Proposition 3.1
in \cite{F1}),
\begin{eqnarray*}
\frac{\partial(\mathbf{w}^{\varepsilon})_k}{\partial x_i}(t,x)
&=&\int_{\mathbb{R}^3}\frac{\partial G^{\nu}_t}{\partial x_i}(x-y)
(w_0)_k(y) \,dy+\int_0^t \int_{\mathbb{R}^3} \frac{\partial G^{\nu
}_t}{\partial x_i}(x-y) \mathbf{g}(0,y)\,dy \,ds\\
&&{} - \int_0^t\sum_{j=1}^3\int_{\mathbb{R}^3}\frac{\partial G^{\nu}_{t-s}}
{\partial x_i}
(x-y) \biggl[\mathbf{K}^{\varepsilon}(\mathbf{w}^{\varepsilon
})_j(s,y)\,\frac{
\partial\mathbf{w}^{\varepsilon}_k(s,y)}{\partial y_j}
\\
&&\hspace*{129.3pt}{} - \mathbf{w}^{\varepsilon}_j(s,y)\, \frac{\partial
\mathbf{K}^{\varepsilon}(\mathbf{w}^{\varepsilon})_k(s,y)}{\partial
y_j} \biggr]\,dy \,ds
\end{eqnarray*}
for $k=1,2,3$. Notice now that, thanks to the estimates
(\ref{unifbound11}), Lemma \ref{aproxims}(ii) also holds
with $p$ replaced by $\tilde{p}$. By estimates as those in the
proof of Theorem \ref{exist1}(i) and using Lemma
\ref{aproxims}(ii) and estimates (\ref{intboundqtilde}) and
(\ref{unifbound11}), we then have
\begin{eqnarray*}
&&\|\nabla\mathbf{w}^{\varepsilon}(t)-\nabla
\mathbf{w}(t)\|_{\tilde{p}}
\\
&&\qquad\leq
C\int_0^t(t-s)^{-{3}/({2\tilde{p}})}s^{-{1/2}}
[\|\mathbf{w}^{\varepsilon}(s)-\mathbf{w}(s)\|_{\tilde{p}}\\
&&\qquad\quad\hspace*{111.6pt}{}+\|\mathbf{K}^{\varepsilon}(\mathbf{w})(s)-\mathbf
{K}(\mathbf{w})(s)\|_{\tilde{q}}
]\,ds\\
&&\qquad\quad{}+C\int_0^t(t-s)^{-{3}/({2\tilde{p}})} [\|\nabla
\mathbf{w}^{\varepsilon}(s)+\nabla\mathbf{w}(s)\|_{\tilde{p}}\\
&&\hspace*{100.8pt}\qquad\quad{}+ \|\nabla\mathbf{K}^{\varepsilon}(\mathbf
{w})(s)-\nabla
\mathbf{K}(\mathbf{w})(s)\|_{\tilde{q}} ]\,ds \\
&&\qquad\leq C \varepsilon t^{1-{3}/{\tilde{p}}}+C \varepsilon
\int_0^t
(t-s)^{-{3}/({2\tilde{p}})}s^{-1+{3}/({2\tilde{p}})-{3}/({2p})}\,ds
\\
&&\qquad\quad{} + C \int_0^t (t-s)^{-{3}/({2\tilde{p}})} \|\nabla
\mathbf{w}^{\varepsilon}(s)-\nabla\mathbf{w}(s)\|_{\tilde{p}} \,ds\\
&&\qquad\leq C \varepsilon t^{-{3}/({2\tilde{p}})}+ C \int_0^t
(t-s)^{-{3}/({2\tilde{p}})} \|\nabla
\mathbf{w}^{\varepsilon}(s)-\nabla\mathbf{w}(s)\|_{\tilde{p}} \,ds.
\end{eqnarray*}
Iterating the latter sufficiently many times (using the identity
quoted in the proof of Theorem \ref{exist1}) (i), we obtain
that
%
%
\begin{eqnarray}\label{hola}
\|\nabla\mathbf{w}^{\varepsilon}(t)-\nabla\mathbf{w}(t)\|_{\tilde
{p}} &\leq& C
\varepsilon\bigl(t^{-{3}/({2\tilde{p}})}+1\bigr)\nonumber\\[-8pt]\\[-8pt]
&&{} + C(T) \int_0^t \|\nabla
\mathbf{w}^{\varepsilon}(s)-\nabla\mathbf{w}(s)\|_{\tilde{p}}
\,ds.\nonumber
\end{eqnarray}
Integrating (\ref{hola}) in time and using Gronwall's lemma,
and then inserting the obtained bound in the right-hand side of
(\ref{hola}), we obtain
%
%
\begin{equation}
\|\nabla\mathbf{w}^{\varepsilon}(t)-\nabla\mathbf{w}(t)\|_{\tilde
{p}}\leq C
\varepsilon\bigl(t^{-{3}/({2\tilde{p}})}+1\bigr),
\end{equation}
and the convergence statement for $\nabla\mathbf{w}^{\varepsilon}$
follows.
\end{pf}
\begin{corolario}\label{convvelfield}
Consider fixed real numbers $\tilde{p}\in(\frac{3}{2},3)$ and
$\alpha\in(0,\frac{1}{2})$. Under the assumptions of Theorem
\ref{convortmet}, there exists a constant $\mathbf{C}$ depending
on $\tilde{p},T, \|w_0\|_p, \tn
\mathbf{g}\tn_{0,p,T},\|w_0\|_1+\|\mathbf{g}\|_{1,T}$ and $\alpha$,
such that for
all $n\in\mathbb{N}$,
\begin{eqnarray*}
&&\sup_{t\in[0,T]}\gamma(t) E \bigl( |
\mathbf{K}^{\varepsilon_n}(\tilde{\mu}^{n,\varepsilon
_n,R})(t,x)-\mathbf{u}
(t,x) | \bigr)\\
&&\qquad\leq
\mathbf{C} \biggl(\frac{(\ln n)^{{1}/{3}}}
{n^{{1}/{2}-\alpha}} +\frac{(\ln n)^{{1}/{3}}
}{\sqrt{n}}+\frac{1}{(\ln n)^{{1}/{9}}} \biggr),
\end{eqnarray*}
where
$\gamma(t)=(t^{{3}/({2\tilde{p}})}+t^{1-{3}/{2}(
{1}/{\tilde{p}}-{1}/{p})})$.
\end{corolario}
\begin{pf} For all $(t,x)\in[0,T]\times\mathbb{R}^3$, it holds
that
%
%
\begin{eqnarray}\label{3estimatesu}\qquad
&&|\mathbf{K}^{\varepsilon_n}(\tilde{\mu}^{n,\varepsilon
_n,R})(t,x)-\mathbf{u}
(t,x) |\nonumber\\
&&\qquad
\leq
\Biggl|\mathbf{K}^{\varepsilon_n}(\tilde{\mu}^{n,\varepsilon_n,R})(t,x)
\nonumber\\
&&\qquad\quad\hspace*{2pt}{}-\frac{1}{n}\sum_{i=1}^n K_{\varepsilon_n}(x-
X_t^{i,\varepsilon_n,R})\wedge
(\chi_R(\Phi_t^{i,\varepsilon_n,R}))h(\tau,X_0^i)\mathbf{1}_{\{
\tau\geq t\}
} \Biggr|\nonumber\\[-8pt]\\[-8pt]
&&\qquad\quad{} + \Biggl|\frac{1}{n}\sum_{i=1}^n K_{\varepsilon_n}(x-
X_t^{i,\varepsilon_n,R})\wedge
(\chi_R(\Phi_t^{i,\varepsilon_n,R}))h(\tau,X_0^i)\mathbf{1}_{\{
\tau\geq t\}
}\nonumber\\
&&\qquad\quad\hspace*{14.4pt}{} -\int_{\mathcal{C}_T} K_{\varepsilon_n}\bigl(x- y(t)\bigr)\wedge
\chi_R(\phi(t))h(\theta,x(0))P^{\varepsilon_n,R}(d\theta,dy,d\phi
) \Biggr|\nonumber\\
&&\qquad\quad{} +|\mathbf{K}^{\varepsilon_n}(\mathbf{w}^{\varepsilon
_n})(t,x)-\mathbf{u}(t,x)|\nonumber
\end{eqnarray}
with $P^{\varepsilon_n,R}$ as in Corollary \ref{convempproc}. By
similar reasons as in (\ref{estimatesf}), the expectation of the
second term is now bounded by
$\frac{1}{\sqrt{n}}2M_{\varepsilon_n}R \bar{h}$. With the estimate
in Theorem \ref{propchaosepsR} we get that
\begin{eqnarray*}
&&
E |\mathbf{K}_{\varepsilon_n}(\tilde{\mu}^{n,\varepsilon
_n,R})(t,x)-\mathbf{u}(t,x) |
\\
&&\qquad\leq (L_{\varepsilon_n} R
+M_{\varepsilon_n})\bar{h}\frac{1}{\sqrt{n}}
C(\varepsilon_n,R,\bar{h},T)
\\
&&\qquad\quad{} +\frac{2M_{\varepsilon_n}R \bar{h} }{\sqrt{n}} +\|\mathbf
{K}^{\varepsilon
_n}(\mathbf{w}^{\varepsilon_n})(t)-\mathbf{K}(\mathbf{w})(t)\|
_{\infty}.
\end{eqnarray*}
Thus, from the estimates for $L_{\varepsilon}$ and
$M_{\varepsilon}$ we deduce that for fixed $\tilde{p}\in
(\frac{3}{2},3)$,
\begin{eqnarray*}
&&
E |\mathbf{K}_{\varepsilon_n}(\tilde{\mu}^{n,\varepsilon
_n,R})(t,x)-\mathbf{u}(t,x) |
\\
&&\qquad\leq C (1+R \bar{h}T)(R\bar{h}T)\frac{(c\ln n)^{{1/3}}}
{n^{{1}/{2}-\alpha}} + C R \bar{h}\frac{(c\ln n)^{{1/3}}
}{\sqrt{n}} \\
&&\qquad\quad{} + \|\mathbf{w}^{\varepsilon_n}(t)-\mathbf{w}(t)\|_{W^{1,\tilde
{p}}}+\|\mathbf{K}
^{\varepsilon_n}(\mathbf{w})(t)-\mathbf{K}(\mathbf{w})(t)\|
_{W^{1,\tilde{q}}},
\end{eqnarray*}
where $\tilde{q}=\frac{3\tilde{p}}{3-\tilde{p}}<\frac{3p}{3-p}$.
We have used here again the Sobolev inclusions quoted in the proof
of Theorem \ref{regularity}, and Lemma \ref{contK}. Now, by Lemmas
\ref{contgradK} and \ref{aproxco1}, one has
\begin{eqnarray*}
\|\nabla\mathbf{K}^{\varepsilon_n}(\mathbf
{w})(t)-\nabla
\mathbf{K}(\mathbf{w})(t)\|_{\tilde{q}} &\leq& C \|
\varphi_{\varepsilon_n}*\mathbf{w}(t)-\mathbf{w}(t)\|_{\tilde
{q}}\leq C
\varepsilon_n \|\nabla\mathbf{w}(t)\|_{\tilde{q}}\\
&\leq& C
t^{-1+{3/2}({1}/{\tilde{p}}-{1}/{p})}\varepsilon_n,
\end{eqnarray*}
where we have also used part (i) of Theorem \ref{regularity}
in the last inequality. On the other hand,
\[
\| \mathbf{K}^{\varepsilon_n}(\mathbf{w})(t)-
\mathbf{K}(\mathbf{w})(t)\|_{\tilde{q}}\leq C \|
{\varphi_{\varepsilon_n}*\mathbf{w}(t)}-\mathbf{w}(t)\|_{\tilde
{p}}\leq C
\varepsilon_n \|\nabla\mathbf{w}(t)\|_{\tilde{p}}\leq C t^{-{1/2}}
\varepsilon_n
\]
thanks to the estimate (\ref{unifbound11}). From
the previous estimates and Lemmas \ref{aproxims} and
\ref{convgradw}, we deduce that
\begin{eqnarray*}
&&E |\mathbf{K}_{\varepsilon_n}(\tilde{\mu}^{n,\varepsilon
_n,R})(t,x)-\mathbf{u}(t,x) |
\\
&&\qquad\leq C \frac{(\ln n)^{{1}/{3}}} {n^{{1}/{2}-\alpha}} + C
\frac{(\ln n)^{{1}/{3}}
}{\sqrt{n}}
\\
&&\qquad\quad{}+C\varepsilon_n\bigl(t^{-{3}/({2\tilde{p}})}+
t^{-{1}/{2}}+t^{-1+{3}/{2}({1}/{\tilde{p}}-{1}/{p})}\bigr),
\end{eqnarray*}
and the statement follows.
\end{pf}

\section{Convergence rate under additional
regularity assumptions}\label{sec6}

Let us finally explain how the convergence rate can be slightly
improved by assuming further regularity of the data $w_0$ and
$\mathbf{g}$. Since it is an adaptation of the developments in the
previous sections, we only sketch the main arguments.

First, it is possible to show that if the data $w_0$ and $\mathbf{g}$ are
such that
%
%
\begin{equation}\label{adreghyp}
\|w_0\|_{W^{m,p}},\qquad \sup_{t\in[0,T]}\|\mathbf{g}(t)\|
_{W^{m,p}}<\infty
\end{equation}
for some integer $m\geq1$, then the mild solutions
$\mathbf{w}^{\varepsilon}$, $\varepsilon\geq0$, given by Theorem
\ref{exist1} belong to the space $\mathbf{F}_{m+1,p,T}$ of functions
$\mathbf{v}(t)$ such that
\[
\sum_{i=1}^{m-1}\tn D^i \mathbf{v}\tn_{0,p,T} + \tn D^{m} \mathbf
{v}\tn
_{1,p,T}<\infty,
\]
where $D^i$ stands for the $i$th order space derivative. To prove
this, one easily first checks that $\mathbf{w}_0$ belongs to that space,
since the successive derivatives in the convolutions the heat
kernel can be applied to the data $w_0$ and $\mathbf{g}$. On the other
hand, on can show by induction that the bilinear operators
$\mathbf{B}^{\varepsilon}$ are continuous in $\mathbf{F}_{m+1,p,T}$,
and more
generally, in the naturally generalized versions
$\mathbf{F}_{m+1,r,(T;p)}$ of the space $\mathbf{F}_{1,r,(T;p)}$.
That is, the
spaces of functions $\mathbf{v}$ such that
\[
\sum_{i=1}^{m-1}\tn D^i \mathbf{v}\tn_{0,r,(T;p)} + \tn D^{m}
\mathbf{v}
\tn_{1,r,(T;p)}
\]
is finite. From this, one gets a local existence
result in the space $\mathbf{F}_{m+1,p,T}$, from which a regularity result
can be obtained by arguments that can be adapted from those in the
proof Theorem 3.2 in \cite{F1}. Moreover, one also checks that
the norms $\tn\mathbf{w}^{\varepsilon}\tn_{m+1,r,(T;p)}$ are bounded
uniformly in $\varepsilon\geq0$.

Now, we impose additional conditions on the regularizing kernel
$\varphi$, namely:
\begin{longlist}
\item$\int_{\mathbb{R}^3}\varphi(x)\,dx=1$.
\item$\int_{\mathbb{R}^3}|x|^{m+1}|\varphi(x)|\,dx<\infty$.
\item$\int_{\mathbb{R}^3}x_{i_1}\cdots x_{i_r} \varphi(x)
\,dx=0$ for
all $i_1,\ldots,i_r\in\{1,2,3\}$ and $r\leq m$.
\end{longlist}
Such function is called a cutoff function of order $m+1$. Then,
one has the following approximation result (see Lemma 4.4 in
\cite{Rav}):
\[
\|\varphi_{\varepsilon}* w - w\|_r\leq C \varepsilon^{m+1} \|
D^{m+1}w\|_r
\]
for all $w\in W^{m+1,r}$. Therefore, without any modification, for
such function $\varphi$, the proofs of Lemmas \ref{aproxims} and
\ref{convgradw} yield the same convergence results but at rate
$\varepsilon^{m+1}$.

By following exactly the same steps as in the previous section, we
finally deduce:
\begin{teorema}
Assume the hypotheses of Theorems \ref{convortmet} and, moreover, that
(\ref{adreghyp}) holds for some integer $m\geq1$ and that
$\varphi$ is a cutoff of order $m+1$. Then, we have for all $i\leq
n$,
\begin{eqnarray*}
&&\mathbb{E} \Bigl[\sup_{t\in[0,T]} \{
|X^{i,\varepsilon_n,R,n}_t-X^i_t|+
|\Phi^{i,\varepsilon_n,R,n}_t-\Phi^i_t| \} \Bigr] \\
&&\qquad\leq
C(T,w_0,\mathbf{g},\alpha) \biggl[\frac{1}{n^{{1}/{2}-\alpha
}(\ln
n)^{{1}/{9}}}+\frac{1}{(\ln n)^{({m+1})/{9}}} \biggr]
\end{eqnarray*}
and
\begin{eqnarray*}
&&\sup_{t\in[0,T],x\in\mathbb{R}^3}\gamma(t) E \bigl( |
\mathbf{K}^{\varepsilon_n}(\tilde{\mu}^{n,\varepsilon
_n,R})(t,x)-\mathbf{u}
(t,x) | \bigr)
\\
&&\qquad\leq\mathbf{C} \biggl(\frac{(\ln n)^{{1}/{3}}}
{n^{{1}/{2}-\alpha}} +\frac{(\ln n)^{{1}/{3}}
}{\sqrt{n}}+\frac{1}{(\ln n)^{({m+1})/{9}}} \biggr),
\end{eqnarray*}
where $\gamma(t)$ was defined in Corollary \ref{convvelfield},
where the
constants now, moreover, depend on $m$.
\end{teorema}

\section*{Acknowledgments}
I would like to thank Mireille Bossy for suggesting me the use of
cutoff techniques in \cite{Rav}. I also thank an anonymous referee for
carefully reading this work, and for helpful suggestions that allowed
me to improve its presentation.

%

%
\printaddresses

\end{document}